\documentclass[12pt]{iopart}

\usepackage{epsfig,verbatim}
\usepackage{amssymb}
\usepackage{amsfonts}

\newcommand{\eps}{\varepsilon}
\newcommand{\dis}{\displaystyle}
\newcommand{\rens}{\mathbb{R}}
\newcommand{\R}{\mathbb{R}}
\newcommand{\E}{\mathbb{E}}
\newcommand{\zens}{\mathbb{Z}}
\newcommand{\torus}{\mathbb{T}}
\newcommand{\dive}{\mbox{div }}

\newcommand{\text}[1]{\mbox{#1} }
\newcommand{\eqref}[1]{(\ref{#1})}
\newcommand{\tb}{\widetilde{b}}
\newcommand{\tsigma}{\widetilde{\sigma}}
\newcommand{\of}{\overline{f}}
\newcommand{\og}{\overline{g}}
\newcommand{\zbar}{\overline{y}}
\newcommand{\xib}{\alpha}

\newtheorem{proposition}{\bf Proposition}[section]
\newtheorem{lemma}{\bf Lemma}[section]
\newtheorem{coro}{\bf Corollary}[section]
\newtheorem{remark}{\bf Remark}[section]
\newtheorem{defi}{\bf Definition}[section]
\newenvironment{proof}{\noindent \textit{Proof:}}{\hfill $\square$}

\begin{document}

\title{Effective dynamics using conditional expectations}
\author{Fr\'ed\'eric Legoll$^{1,2}$ and Tony Leli\`evre$^{3,2}$}
\address{
$^1$
Universit\'e Paris-Est,
Institut Navier, LAMI, \'Ecole des Ponts,
6 et 8 avenue Blaise Pascal, 77455 Marne-La-Vall\'ee Cedex 2, France
}
\address{
$^2$
INRIA Rocquencourt, MICMAC Team-Project, Domaine de Voluceau, B.P. 105,
78153 Le Chesnay Cedex, France
}
\address{
$^3$
Universit\'e Paris-Est,
CERMICS, \'Ecole des Ponts,
6 et 8 avenue Blaise Pascal, 77455 Marne-La-Vall\'ee Cedex 2, France
}
\eads{\mailto{legoll@lami.enpc.fr},
\mailto{lelievre@cermics.enpc.fr}}

\begin{abstract}
The question of coarse-graining is ubiquitous in molecular dynamics. In
this article, we are interested in deriving effective properties for the
{\em dynamics} of a coarse-grained variable $\xi(x)$,
where $x$ describes the configuration of the system in a high-dimensional space $\R^n$,
and $\xi$ is a smooth function with value in $\R$ (typically a reaction coordinate). It
is well known that, given a Boltzmann-Gibbs distribution on $x \in \R^n$, the equilibrium
properties on $\xi(x)$ are completely determined by the free energy. On the other hand,
the question of the effective dynamics on $\xi(x)$ is much more difficult to address.
Starting from an overdamped Langevin equation on $x \in \R^n$, we propose an effective dynamics
for $\xi(x) \in \R$ using conditional expectations. Using entropy methods, we give sufficient
conditions for the time marginals of the effective dynamics to be close to the original
ones. We check numerically on some toy examples that these sufficient conditions 
yield an effective dynamics which accurately reproduces the
residence times in the potential energy wells.  
We also discuss the accuracy of the effective dynamics in a pathwise sense, and the
relevance of the free energy to build a coarse-grained dynamics.
\end{abstract}


\maketitle

\section{Motivation}

In molecular dynamics, two types of quantities are typically of interest: averages with respect to the canonical ensemble (thermodynamic quantities, such as stress, or heat capacity), and averages of functionals over paths (dynamic quantities, like viscosity, diffusion coefficients or rate constants). In both cases, the question of coarse-graining is relevant, in the sense that the considered functionals typically depend only on a few variables of the system (collective variables, or reaction coordinates) so that it would be interesting to obtain coarse-grained models on these variables.

\subsection{Coarse-graining of thermodynamic quantities}

Computing canonical averages is a standard task in molecular
dynamics. For a molecular system whose atom positions are described by a
vector $x \in \rens^n$, these quantities read 
\begin{equation}
\label{eq:average}
\int_{\R^n} \Phi(x) \, d\mu
\end{equation}
where $\Phi : \R^n \to \R$ is the observable of interest and $\mu$ is the
Boltzmann-Gibbs measure,
\begin{equation}
\label{eq:mu}
d\mu = Z^{-1} \exp(-\beta V(x)) \, dx,
\end{equation}
where $V$ is the potential energy of the system, $\beta$ is
proportional to the inverse of the system temperature, and $\dis{
Z=\int_{\R^n} \exp(-\beta V (x)) \, dx
}$
is a normalizing constant. Typically, $x$ represents the position of $N$
three-dimensional particles, hence $x \in \rens^n$ with $n = 3N$. All
the results we prove 
are also satisfied if $x \in \torus^n$, where $\torus = \R / \zens$
denotes the one-dimensional torus. 

As mentioned above, observables of interest are often function of only part of the
variable $x$. For example, $x$ denotes the positions of {\em all}
the atoms of a protein and of the solvent molecules around, and the
quantity of interest is only a particular angle between some atoms in
the protein, because this angle characterizes the conformation of the
protein (and thus the potential energy well in which the system is is
completely determined by the knowledge of this quantity of interest). We
thus introduce the so-called {\em reaction coordinate} 
$$
\xi: \R^n \to \R,
$$
which contains all the information we are interested in~\footnote{In this
article, we do not address the difficult question of how to find a good
reaction coordinate. See for instance~\cite{michalak} for some
discussion on that point.}.
Throughout this article, we assume that 
$$
  \text{{\bf [H1]}~~~}
\begin{array}{c}
\text{$\xi$ is a smooth {\em scalar} function such that,}
\\
\text{for all } x \in \R^n, \ 0 < m \le |\nabla \xi(x)| \le M < \infty.
\end{array}
$$
We have supposed that $\xi$ is a scalar function. 
It is not clear to us whether the results of this article can be
generalized to the case of a multi-dimensional reaction coordinate. 

To this function $\xi$ is naturally associated an effective
energy $A$, called the {\em free energy}, such that
\begin{equation}
\label{eq:xi}
d(\xi \star \mu) = \exp(-\beta A(z))\, dz,
\end{equation}
where $\xi \star \mu$ denotes the image of the measure $\mu$ by $\xi$. In
other words, for any test function $\Phi : \R \to \R$,
\begin{equation}
\label{eq:stat}
\int_{\R^n} \Phi(\xi(x)) \, Z^{-1} \exp(-\beta V(x)) \, dx 
= 
\int_\R \Phi(z) \, \exp(-\beta A(z))\, dz.
\end{equation}
Expressions of $A$ and its derivative are given below (see Section
\ref{sec:notation}).

The interpretation of \eqref{eq:stat} is that, when $X$ is distributed
according to the Boltzmann measure \eqref{eq:mu}, then $\xi(X)$ is
distributed according to the measure $\exp(-\beta A(z))\, dz$. 
Hence, the free energy $A$ is a relevant quantity for computing thermodynamic quantities, namely canonical averages.

In conclusion, the question of coarse-graining thermodynamic quantities amounts to computing the free energy, and there are several efficient methods to perform such calculations (see for example~\cite{chipot-pohorille-07}). There are also interesting questions related to computing approximations of the free energy, especially when the number of reaction coordinates is large, for example in polymer science, but this is not the subject of this article.

\subsection{Coarse-graining of dynamical quantities}

The objective of this work is to address some issues related to the {\em
dynamics} of the system, and how to coarse-grain it. In short, we
aim at designing a dynamics that approximates the path $t \mapsto
\xi(X_t)$, where $\xi$ is the above reaction coordinate.

To make this question precise, we first have to {\em choose} the full dynamics, which will be the reference one.
In the following, we consider the overdamped Langevin dynamics on state space~$\R^n$ (we will discuss this choice below),
\begin{equation}
\label{eq:X}
dX_t = - \nabla V(X_t) \, dt + \sqrt{2 \beta^{-1}} \, d W_t,
\quad X_{t=0} = X_0,
\end{equation}
where $W_t$ is a standard $n$-dimensional Brownian motion.
Under suitable assumptions on~$V$, this dynamics is ergodic with respect
to the Boltzmann-Gibbs measure \eqref{eq:mu}. Hence, for $\mu$-almost
all initial conditions $X_0$,
\begin{equation}
\label{eq:ergo}
\lim_{T \to \infty} \frac{1}{T} \int_0^T \Phi(X_t) \, dt = 
\int_{\R^n} \Phi(x) \, d\mu
\end{equation}
almost surely. In practice, this convergence is very slow, due to some
metastabilities in the dynamics: $X_t$ samples a given well of the
potential energy for a long time, before hoping to some other well of
$V$. 

An important dynamical quantity we will consider below is the
average residence time, that is the mean time that the system spends in a
given well, before hoping to another one, when it follows the dynamics
\eqref{eq:X}. Typically, the wells are fully described through $\xi$
($x$ is in a given well if and only if $\xi(x)$ is in a given interval),
so that these times can be obtained from the knowledge of the time
evolution $\xi(X_t)$, which is expensive to compute since it means 
  simulating the full system.

In this article, our aim is
  twofold. First, we would like to propose a 
one-dimensional dynamics of the form
\begin{equation}
\label{eq:eff_dyn}
dy_t = b(y_t) \, dt + \sqrt{2 \beta^{-1}} \, \sigma(y_t) \, d B_t,
\end{equation}
where $B_t$ is a standard one-dimensional Brownian motion and $b$ and
$\sigma$ are scalar functions, such that $\left( y_t \right)_{0 \leq t \leq T}$
is a good approximation (in a sense to be made precise below) of
$\left( \xi(X_t)\right)_{0 \leq t \leq T}$. Hence, the dynamics \eqref{eq:eff_dyn} can be 
thought of as a coarse-grained, or {\em effective}, dynamics for the quantity
of interest. A natural requirement is that \eqref{eq:eff_dyn} preserves
equilibrium quantities, {\em i.e.} it is ergodic with respect to
$\exp(-\beta A(z))\, dz$, the equilibrium measure of $\xi(X)$, but we
typically ask for more than that. For example, we would like to be able to
recover residence times in the wells from~\eqref{eq:eff_dyn}, hence
bypassing the expensive simulation of $\xi(X_t)$ (see Section~\ref{sec:num_D2} for
some numerical results on that quantity).  

Second, we would like to investigate the relation between
\eqref{eq:eff_dyn} and the coarse-grained dynamics 
\begin{equation}
\label{eq:zbar}
d\zbar_t = - A'(\zbar_t) \, dt + \sqrt{2 \beta^{-1}} \, d B_t,
\end{equation}
which is indeed a one-dimensional dynamics, driven by the free energy,
and ergodic for $\exp(-\beta A(z))\, dz$. In other words, what is the
dynamical content of the free energy? 
This second question stems from the fact that practitioners
often look at the free energy profile ({\em i.e.} the function $z \mapsto
A(z)$) to get an idea of the dynamics of transition (typically the
transition time) between one region indexed by the reaction coordinate
(say for example $\left\{ x \in \R^n; \ \xi(x) \le z_0 \right\}$) and
another one (for example $\left\{ x \in \R^n; \ \xi(x) > z_0
\right\}$). If $\xi(X_t)$ follows a dynamics which is close to
\eqref{eq:zbar}, then the Transition State Theory says that residence
times are a function of the free energy
barriers~\cite{kramers,kramers_review}, and then it makes sense to look
at the free energy to compute some dynamical properties.
It is thus often assumed that there is some dynamical information in 
the free energy $A$.

\medskip

The difficulty of the question we address stems from the fact
that, in general, $t \to \xi(X_t)$ is not a Markov process: this is a closure problem. A first
possibility is to try and approximate $\xi(X_t)$ by a process which has
some memory in time, typically a generalized Langevin equation (see for
instance~\cite{darve09,FU_gle}, and also~\cite{hartmann_phd}). 
A standard framework is then the Mori-Zwanzig
projection formalism, which is described in details
in~\cite{gks}. 
Note also that, since we are interested in reproducing only some output
function of $X_t$ (namely $\xi(X_t)$), tools from the control theory may
be used. Such an idea has been followed in~\cite{hart_schu_ode,hart_schu_sde}.  

If a time-scale separation is present in the
system, then memory effects may be neglected. In the sequel, we 
make such time-scale separation assumptions (see assumptions {\bf [H2]} and {\bf [H3]} of
Proposition~\ref{prop:D3}), which allow us to approximate $\xi(X_t)$ by
a Markov process of the type~\eqref{eq:eff_dyn}. We use the
framework of logarithmic Sobolev inequalities to write these
assumptions. It has the advantage that we do not assume to {\em a
  priori} know how to split $x$ between fast and slow modes, or to
split the potential energy $V$ between fast and slow terms (otherwise
stated, the time scale separation is encoded in the constants entering
the logarithmic Sobolev inequalities, and not inserted {\em a priori} in
the model). In addition, within this framework, we can handle reaction
coordinates that are nonlinear functions of $x$, the natural cartesian
coordinates of the system (see the numerical simulations reported in
Section~\ref{sec:num_D2}). 

Another possibility is to start from a dynamics which includes an
{\em explicit} small parameter, representing a time scale separation. One may
then apply an averaging principle (see~\cite{hartmann_phd} and the references
therein for more details along this idea; see also~\cite{pav_stu} for a
comprehensive review of the averaging principle, when applied to
deterministic and stochastic differential equations). In
Section~\ref{sec:para}, we consider such a case of potential
energy being the sum of two terms of different stiffness, 
as an {\em example} of
application of our general result (see the potential
energy~\eqref{eq:ex_V}). Note that, even if we explicitly insert a small parameter in
$V$, our model differs from the one considered in~\cite{eve03},
where a small parameter appears in the potential energy {\em and} in the
diffusion coefficient. 

Other strategies are to try and identify fast and slow modes in the
dynamics (see e.g.~\cite{iza08,iza09}), or to postulate a parametric
form for the effective dynamics and to identify its coefficients by
numerical simulation on the complete system~\cite{pokern,yang}. 

We finish this section by a discussion of the choice of the full dynamics. We chose the overdamped Langevin dynamics~\eqref{eq:X}. Other
choices can be made, in particular the Langevin dynamics, which is closer to a Hamiltonian dynamics and can also be seen as a method to sample the canonical measure (see~\cite{comparisonNVT} for a review
of sampling methods of the canonical ensemble, along with a
theoretical and numerical comparison of their performances for molecular
dynamics). From the analysis standpoint, the dynamics we chose is much
simpler, since the diffusion is non-degenerate (in contrast to the Langevin
dynamics, which is an hypoelliptic equation). We do not know whether the theoretical
results presented in this article (such as Proposition~\ref{prop:D3})
can be generalized to the case of the Langevin dynamics. From a
practical viewpoint, it may be possible to use the same strategy starting from
the Langevin dynamics to write another low-dimensional dynamics. We have
not pursued in this direction. As an alternative to continuous time
processes, one can also model the dynamics of a molecular system by a
discrete time Markov chain, for instance in a discrete state space,
where each state represents a different metastable configuration of the
system~\cite{schuette_jcp,schuette_handbook}. 
In that setting, the question of estimating the accuracy of a
coarse-grained dynamics has been addressed in~\cite{schuette_ima}, where
similar bounds as those derived in this article are obtained.

\subsection{Statement of the main results and outline}

We propose a way to derive an effective dynamics of the
form~\eqref{eq:eff_dyn}. This defines a process $(y_t)_{t\ge 0}$, which we compare
with $(\xi(X_t))_{t\ge 0}$, where $X_t$ satisfies~\eqref{eq:X}. Three
quantities can be typically considered to estimate the distance 
between $y_t$ and $\xi(X_t)$ (on the time interval $[0,T]$):
\begin{itemize}
\item ~[D1] pathwise convergence: $ \E \left(\sup_{t \in (0,T)}
    |\xi(X_t) - y_t|^2 \right)$,
\item ~[D2] convergence of the laws of paths: 
$\|{\mathcal L}(\xi(X_t)_{0 \le t \le T}) - {\mathcal L}((y_t)_{0 \le t
  \le T}) \|_{TV}$, 
\item ~[D3] convergence of time marginals: $\sup_{t \in (0,T)}
  \|{\mathcal L}(\xi(X_t)) - {\mathcal L}(y_t) \|_{TV}$.
\end{itemize}
In the above estimators, 
we have arbitrarily chosen to measure distances between probability
measures by the total variation (TV) norm, but other choices could be
made. Recall that the total variation of a signed measure $\nu$ is
defined by 
$
\| \nu \|_{TV} = \sup_{f \in L^\infty, \| f \|_{L^\infty} \leq 1} 
\int f d\nu.
$
If $\nu$ is a measure on $\R^n$ which has a density with respect to the
Lebesgue measure, then its total variation is just the $L^1$ norm of its
density.

It is clear that a bound in the sense of [D1] implies a bound in the
sense of [D2], which implies a bound in the sense of [D3]. Conversely, 
by the Skorohod theorem, a bound in the sense of
[D2] implies a bound in the sense of [D1], for some well chosen
realizations of $W_t$ and $B_t$ (the brownian motions in \eqref{eq:X}
and \eqref{eq:eff_dyn}), but this theorem is not
constructive. 
The most relevant criterion in practice is [D2]. Indeed, the criterion
[D3] does not account for the correlations in time of the process, which are important
to understand its dynamical properties. On the other hand, the
pathwise convergence criterion [D1] is too strong: practionners in
molecular dynamics are rarely interested in the trajectory {\em per
  se}. Moreover, [D2] implies the convergence of the law of escape times
(hence of residence times in the wells), at least if the escape time is
(almost surely) a continuous function of paths, which holds under some
regularity assumptions (see~\cite[Exercise 3.9.10]{bichteler-02}).

Our first objective is to propose, in a general case, some sufficient
conditions on the reaction coordinate for a bound of type [D3] to be
satisfied. We are actually able to derive an estimate of the difference between the time marginals which is {\em uniform in time}. Next, on a toy-model, we investigate, both 
theoretically and numerically, if these conditions are sufficient and
necessary for [D1] and [D2] to hold. 

The article is organized as follows. In Section~\ref{sec:eff_dyn}, after
introducing some notation and recalling some basic relations concerning
the free energy, we propose a natural coarse-graining procedure, which
enables us to obtain an effective dynamics of type
\eqref{eq:eff_dyn}, where the functions $b$ and $\sigma$ can easily be
computed (see Equations~\eqref{eq:b},~\eqref{eq:sigma}
and~\eqref{eq:y}). In Section \ref{sec:preuve}, we prove that the
solution $y_t$ of the effective dynamics~\eqref{eq:y} is indeed a good
approximation of $\xi(X_t)$, 
in the sense [D3]. Our argument relies on entropy techniques, and is
very much inspired by \cite{otto09,dizdar}. 
In Section~\ref{sec:num_D2}, we present some numerical results obtained
on a simple model, where we compare residence times in the potential
energy wells as predicted by the reference dynamics~\eqref{eq:X} and by
the one-dimensional reduced dynamics~\eqref{eq:y}. Section~\ref{sec:path} is
dedicated to establishing error estimates in the sense [D1] of pathwise
convergence, in a specific case. These estimates are illustrated by numerical
simulations.  

\section{A ``natural'' coarse-graining procedure}
\label{sec:eff_dyn} 

\subsection{Notation}
\label{sec:notation}

We gather here some useful notation and results. Let $\Sigma_z$ be
the submanifold of $\R^n$ of positions at a fixed value of the reaction coordinate:
$$
\Sigma_z= \{ x \in \R^n; \, \xi(x) =  z \}.
$$
Let us introduce $\mu_{\Sigma_z}$, which is the probability measure
$\mu$ conditioned at a fixed value of the reaction coordinate:
\begin{equation}
\label{eq:mu_z}
d\mu_{\Sigma_z} = 
\frac{\exp(-\beta V) \, |\nabla \xi|^{-1} \, d\sigma_{\Sigma_z}}
{\dis{ \int_{\Sigma_z} \exp(-\beta V) \, |\nabla \xi|^{-1} \,
    d\sigma_{\Sigma_z}}}, 
\end{equation}
where the measure $\sigma_{\Sigma_z}$ is the Lebesgue measure
on~$\Sigma_z$ induced by the Lebesgue measure in the ambient Euclidean
space $\R^n \supset \Sigma_z$. 

We recall the following expressions for the free energy $A$ and its
derivative $A'$, also called the {\em mean force} (see~\cite{c_l_eve}):
\begin{equation}
\label{eq:A}
A(z) = - \beta^{-1} \ln \left( \int_{\Sigma_z} Z^{-1} \, \exp(-\beta V)
  \, |\nabla \xi|^{-1} \, d\sigma_{\Sigma_z} \right)
\end{equation}
and
\begin{equation}
\label{eq:A'}
A'(z) = \int_{\Sigma_z} F \, d\mu_{\Sigma_z},
\end{equation}
where $F$ is the so-called {\em local mean force}:
\begin{equation}
\label{eq:F}
F=\frac{\nabla V \cdot \nabla \xi}{|\nabla \xi|^2}  - \beta^{-1} \, \dive
\left( \frac{\nabla \xi}{|\nabla \xi|^2} \right).
\end{equation}
In view of~\eqref{eq:A}, note that~\eqref{eq:mu_z} reads 
\begin{equation}
\label{eq:mu_z_bis}
d\mu_{\Sigma_z} = 
\frac{\exp(-\beta V) \, |\nabla \xi|^{-1} \, d\sigma_{\Sigma_z}}
{Z \exp(-\beta A)}.
\end{equation}
These expressions can be obtained by the co-area
formula~\cite{evans-gariepy-92}, which we now recall:
\begin{lemma}
For any smooth function $\Phi: \R^{n}  \to \R$,
\begin{equation}
\label{eq:co-area}
\int_{\R^{n}} \Phi(x) \, |\nabla \xi(x)| \, dx = \int_\R
\int_{\Sigma_z} \Phi \, d\sigma_{\Sigma_z} \, dz.
\end{equation}
\end{lemma}
\begin{remark}[Co-area formula and conditioning]
\label{rem:coarea}
The co-area formula shows that if the random variable $X$ has law
$\psi(x)\, dx$ in~$\R^n$, then $\xi(X)$ has law $\psi^\xi(z) \,dz$, with
$$
\psi^\xi(z)=\int_{\Sigma_z} \psi \, |\nabla \xi|^{-1} \,
d\sigma_{\Sigma_z}.
$$
It also shows that the law of $X$ conditioned to a fixed value $z$ of
$\xi(X)$ is $\mu_{\Sigma_z}$, where $\mu_{\Sigma_z}$ is 
defined by~\eqref{eq:mu_z}.
The measure $|\nabla \xi|^{-1} d\sigma_{\Sigma_z}$ is sometimes denoted
by $\delta_{\xi(x)-z}$ in the literature. 
\end{remark}

From the co-area formula, we get the following result:
\begin{lemma}
\label{lem:chi_bar_prime}
For any smooth function $\chi:\R^n \to \R$, consider
$$
\chi^\xi(z)=\int_{\Sigma_z} \chi \, |\nabla \xi|^{-1} \,
d\sigma_{\Sigma_z}.
$$
The derivative of $\chi^\xi$ reads: 
$$
\frac{d\chi^\xi}{dz}(z) = \int_{\Sigma_z} \left[
\frac{\nabla \xi \cdot \nabla \chi}{| \nabla \xi |^2} 
+ 
\chi \, \mbox{\rm div } \left( \frac{\nabla \xi}{| \nabla \xi |^2} \right) 
\right] \,
|\nabla \xi|^{-1} \, d\sigma_{\Sigma_z}.
$$
\end{lemma}
\begin{proof}
For any smooth test function $g:\R \to \R$, we obtain, using the co-area
formula~(\ref{eq:co-area}), that
$$
\int_\R \chi^\xi(z) \, g'(z) \, dz
=
\int_\R \int_{\Sigma_z} \chi \, |\nabla \xi|^{-1} \, g'(z) \, 
d\sigma_{\Sigma_z} \, dz
=
\int_{\R^n} \chi(x) \, g'(\xi(x)) \, dx.
$$
Hence, 
\begin{eqnarray*}
\dis
\int_\R \chi^\xi(z) \, g'(z) \, dz
&=& \dis
\int_{\R^n} \chi(x) \, g'(\xi(x)) \, dx
\\
&=& \dis
\int_{\R^n} \chi \, |\nabla \xi|^{-2} \ \nabla \xi
\cdot \nabla (g \circ \xi)
\\
&=& \dis
- \int_{\R^n} g \circ \xi \ \dive \left( \chi \, |\nabla \xi|^{-2} \, \nabla \xi \right) 
\\
&=& \dis
-\int_{\R} g(z) \int_{\Sigma_z} \dive \left( \chi \, |\nabla \xi|^{-2} \, \nabla \xi \right)
\, \frac{d\sigma_{\Sigma_z}}{|\nabla \xi|} \ dz,
\\
&=& \dis
-\int_{\R} g(z) \int_{\Sigma_z} 
\left[ \frac{\nabla \xi \cdot \nabla \chi}{|\nabla \xi|^2} 
+ \chi \, \dive \left( \frac{\nabla \xi}{|\nabla \xi|^2} \right)
\right]
\, \frac{d\sigma_{\Sigma_z}}{|\nabla \xi|} \ dz,
\end{eqnarray*}
which yields the result.
\end{proof} 

\subsection{A non-closed equation}

Consider $X_t$ that solves~\eqref{eq:X}. 
By a simple It\^o computation, we have
\begin{equation}
\label{eq:xi(X)}
d \xi(X_t) = \left(- \nabla V \cdot \nabla \xi + \beta^{-1} \Delta \xi
\right) (X_t) \, dt + \sqrt{2 \beta^{-1}} \ |\nabla \xi|(X_t) \, dB_t
\end{equation}
where $B_t$ is the one-dimensional Brownian motion
\begin{equation}
\label{eq:dB}
dB_t = \frac{\nabla\xi}{|\nabla\xi|}(X_t) \cdot dW_t.
\end{equation}
Of course, equation~\eqref{eq:xi(X)} is not closed. Following
Gy\"ongy~\cite{gyongy-86}, a simple closing 
procedure is to consider $\widetilde{y}_t$ solution to
\begin{equation}
\label{eq:ty}
d \widetilde{y}_t = \tb(t, \widetilde{y}_t) \, dt + \sqrt{2 \beta^{-1}} \
\tsigma(t, \widetilde{y}_t) \, dB_t,
\end{equation}
where
\begin{equation}
\label{eq:tb}
\tb(t,y)=\E\left[ \left( - \nabla V \cdot \nabla \xi + \beta^{-1} \Delta \xi
\right) (X_t) \ | \ \xi(X_t)=y \right]
\end{equation}
and
\begin{equation}
\label{eq:tsigma}
\tsigma^2(t,y) = \E\left[|\nabla \xi|^2(X_t) \ | \ \xi(X_t)=y \right].
\end{equation}
Note that $\tb$ and $\tsigma$ depend on $t$, since these are expected values conditioned on the fact that $\xi(X_t)=y$, where the
probability distribution function of $X_t$ of course depends on $t$. 

As shown in~\cite{gyongy-86}, 
this procedure is exact from the point of view of time marginals,
{\em i.e.} [D3] in our above classification. This is stated in the following lemma:

\begin{lemma}
\label{lem:gyongy}
The probability distribution function $\psi^\xi$ of $\xi(X_t)$, where $X_t$ satisfies~\eqref{eq:X}, satisfies
the Fokker-Planck equation associated to~\eqref{eq:ty}: 
\begin{equation}
\label{eq:tFP}
\partial_t \psi^\xi = \partial_z\left( - \tb \ \psi^\xi + \beta^{-1}
  \partial_z (\tsigma^2 \psi^\xi )\right).
\end{equation}
\end{lemma}

\begin{proof}
Let us denote $\psi(t,x)$ the probability distribution function of
$X_t$. It satisfies the Fokker-Planck equation
\begin{equation}
\label{eq:FP_full}
\partial_t \psi = \dive \left( \nabla V \, \psi + \beta^{-1} \nabla \psi \right).
\end{equation}
In view of Remark~\ref{rem:coarea}, the probability distribution
function $\psi^\xi(t,z)$ of $\xi(X_t)$ is given by 
$$
\psi^\xi(t,z)=\int_{\Sigma_z} \psi(t,\cdot) \, |\nabla \xi|^{-1} \,
d\sigma_{\Sigma_z}.
$$
Using Lemma~\ref{lem:chi_bar_prime} with $\chi \equiv \psi(t,\cdot)$, we obtain 
\begin{equation}
\label{eq:derpsi}
\partial_z \psi^\xi(t,z)=
\int_{\Sigma_z} 
\left( 
\frac{\nabla \xi \cdot \nabla \psi(t,\cdot)}{|\nabla \xi|^2} 
+ \dive \left(\frac{\nabla \xi}{|\nabla \xi|^2} \right) \psi(t,\cdot) 
\right)
|\nabla \xi|^{-1} \ d\sigma_{\Sigma_z}.
\end{equation}
By definition, we have the following expressions for $\tb$ and
$\tsigma$ in terms of $\psi$:
\begin{eqnarray*}
\tb(t,z) &=& \frac{1}{\psi^\xi(t,z)} 
\int_{\Sigma_z} \left(- \nabla V \cdot \nabla \xi + \beta^{-1} \Delta \xi
\right) |\nabla \xi|^{-1} \ \psi \ d \sigma_{\Sigma_z},
\\
\tsigma^2(t,z) &=& \frac{1}{\psi^\xi(t,z)} 
\int_{\Sigma_z} |\nabla \xi| \ \psi \ d \sigma_{\Sigma_z}.
\end{eqnarray*}
Using again Lemma~\ref{lem:chi_bar_prime} with $\chi \equiv |\nabla
\xi|^2 \ \psi(t,\cdot)$, we obtain 
\begin{equation}
\label{eq:dersigpsi}
\partial_z (\tsigma^2 \, \psi^\xi)
=
\partial_z \int_{\Sigma_z} |\nabla \xi| \, \psi \, d \sigma_{\Sigma_z}
=
\int_{\Sigma_z} 
\left( \nabla \xi \cdot \nabla \psi + \psi \, \Delta \xi \right)
|\nabla \xi|^{-1} \ d\sigma_{\Sigma_z}.
\end{equation}
Let us now prove a variational formulation of \eqref{eq:tFP}. For any
test function $g$, we have
\begin{equation*}
\hspace{-2.5cm}
\begin{array}{rcl}
\dis
\frac{d}{dt} \int_{\R} \psi^\xi(t,z) \, g(z) \, dz
&=& \dis
\frac{d}{dt} \int_{\R^n} \psi(t,x) \, g(\xi(x)) \, dx
\\
&=& \dis
\int_{\R^n} \dive\left( \psi \nabla V + \beta^{-1} \nabla \psi
\right) g \circ \xi \ dx
\\
&=& \dis
- \int_{\R^n} \left( \psi \nabla V + \beta^{-1}
  \nabla \psi \right) \cdot \nabla \xi \ \ g' \circ \xi \ dx
\\
&=& \dis
- \int_{\R}  \int_{\Sigma_z} |\nabla \xi|^{-1} 
\left( \psi \nabla V \cdot \nabla \xi + \beta^{-1} \nabla \psi \cdot
  \nabla \xi \right) \, d \sigma_{\Sigma_z} \ g'(z) \, dz
\\
&=& \dis
- \beta^{-1} \int_{\R} \partial_z (\tsigma^2 \psi^\xi) \, g'(z) \,
dz 
\\ 
&& \dis
+ \int_{\R} \int_{\Sigma_z} |\nabla \xi|^{-1} 
\left(- \nabla V \cdot \nabla \xi + \beta^{-1} \Delta \xi \right) \,
\psi \, d \sigma_{\Sigma_z} \, g'(z) \, dz
\\
&=& \dis
- \beta^{-1} \int_{\R} \partial_z (\tsigma^2 \, \psi^\xi) \, g'(z)  \,
dz + \int_{\R} \tb \ \psi^\xi \ g'(z) \, dz.
\end{array}
\end{equation*}
This shows that $\psi^\xi$ satisfies~\eqref{eq:tFP}.
\end{proof}

\subsection{A closed effective dynamics}
\label{sec:closed}

The problem with equation~\eqref{eq:ty} is that the functions $\tb$
and $\tsigma$ are very complicated to compute, since they involve
the full knowledge of $\psi$. Therefore, one cannot
consider~\eqref{eq:ty} as a reasonable closure. A natural
simplification is to
consider a time-independent approximation of the functions $\tb$
and $\tsigma$. Considering~\eqref{eq:tb} and~\eqref{eq:tsigma}, we
introduce ($\E_{\mu}$ denoting a mean with respect to the measure $\mu$)
\begin{eqnarray}
b(z)
&=&\E_{\mu} \left[ \left(- \nabla V \cdot \nabla \xi + \beta^{-1} \Delta \xi
\right) (X) \ | \ \xi(X)=z \right],\nonumber \\
&=& \int_{\Sigma_z} \left(- \nabla V \cdot \nabla \xi + \beta^{-1} \Delta \xi
\right) d \mu_{\Sigma_z},\label{eq:b}
\end{eqnarray}
and
\begin{eqnarray}
\sigma^2(z)
&=&\E_{\mu}\left(|\nabla \xi|^2(X) \ | \ \xi(X)=z \right),\nonumber \\
&=&\int_{\Sigma_z} |\nabla \xi|^2 \ d \mu_{\Sigma_z},\label{eq:sigma}
\end{eqnarray}
where $\mu_{\Sigma_z}$ is defined by \eqref{eq:mu_z}. This
simplification especially makes sense if $\xi(X_t)$ is a slow variable,
that is if the characteristic evolution time of $\xi(X_t)$ is much
larger than the characteristic time needed by $X_t$ to sample the
manifold $\Sigma_z$. This is quantified in the sequel. 

In the spirit of \eqref{eq:ty}, we
next introduce the coarse-grained dynamics 
\begin{equation}
\label{eq:y}
\fbox{$\dis{
d y_t = b(y_t) \, dt + \sqrt{2 \beta^{-1}} \ \sigma(y_t) \, dB_t, 
\quad y_{t=0}=\xi(X_0).
}$}
\end{equation}
The Fokker-Planck equation associated to the above dynamics will be
useful. It reads
\begin{equation}
\label{eq:FP}
\partial_t \phi = \partial_z \left( - b \, \phi + \beta^{-1}
  \partial_z (\sigma^2 \, \phi) \right).
\end{equation}

Let us first prove that the dynamics~\eqref{eq:y} is ergodic for the
equilibrium 
measure $\xi \star \mu$. The distance between $y_t$ and
$\xi(X_t)$ is estimated in Section \ref{sec:preuve}. 

\medskip

In view of assumption {\bf [H1]} and of~\eqref{eq:sigma}, we
observe that the diffusion coefficient of~\eqref{eq:y} satisfies
$\sigma(y) \geq m > 0$ for any $y$. Hence, the process defined by~\eqref{eq:y} is irreducible,
and admits a unique invariant probability measure. In the following
lemma, we prove that $\exp(-\beta A(z)) \, dz$ is a stationary measure
for~\eqref{eq:y}. Hence, the process $y_t$ defined by~\eqref{eq:y} is
ergodic with respect to this probability (see Has'minskii~\cite{hash},
Kliemann~\cite{kliemann} and the references therein). 

\begin{lemma}
\label{lem:CG_sampling}
The measure $\xi \star \mu$ on $\R$, which has the density
$\exp(-\beta A)$, is a stationary measure for~\eqref{eq:y}.
\end{lemma}

\begin{proof}
We infer from~\eqref{eq:sigma} and~\eqref{eq:mu_z_bis} that
$$
\sigma^2 \, \exp(-\beta A) = Z^{-1} \int_{\Sigma_z} 
\left| \nabla \xi \right| \, \exp(-\beta V) \, d \sigma_{\Sigma_z}.
$$
Using Lemma~\ref{lem:chi_bar_prime} with 
$\chi \equiv Z^{-1} \left| \nabla \xi \right|^2 \, \exp(-\beta V)$, we
obtain  
\begin{eqnarray}
\beta^{-1} \partial_z( \sigma^2 \, \exp(-\beta A))
\nonumber
\\
\quad
=
\beta^{-1} Z^{-1} \int_{\Sigma_z} \left[ \nabla \xi \cdot
  \nabla(\exp(-\beta V)) + \exp(-\beta V) \Delta \xi \right] \ 
|\nabla \xi|^{-1} \ d \sigma_{\Sigma_z},
\nonumber
\\
\quad
=
Z^{-1} \int_{\Sigma_z} \left(- \nabla \xi \cdot
  \nabla V + \beta^{-1} \Delta \xi \right) \exp(-\beta V) \ 
|\nabla \xi|^{-1} \ d \sigma_{\Sigma_z},
\nonumber
\\
\quad
= b \ \exp(-\beta A).
\label{eq:detailed}
\end{eqnarray}
As a consequence of the above equation, \eqref{eq:FP} can be
recast as
\begin{eqnarray}
\label{eq:FP'}
\partial_t \phi &=&
\partial_z \left( - b \, \phi + \beta^{-1}
\partial_z \left(\sigma^2 \, \exp(-\beta A) \, \exp(\beta A) \, \phi \right)
\right)
\nonumber
\\
&=&
\beta^{-1} \partial_z \left[ \sigma^2 \
  \partial_z (\phi \exp(\beta A) ) \ \exp(-\beta A) \right].
\end{eqnarray}
It is now clear that $\phi = \exp(-\beta A)$ is a
stationary solution of the above equation.
\end{proof}

\smallskip

In view of~\eqref{eq:FP'}, we observe that $\phi = \exp(-\beta A)$
is not only a stationary measure for~\eqref{eq:y}, but also satisfies a
detailed balance condition ($(y_t)$ is a reversible process with respect to $\exp(-\beta A(z)) \, dz$). 

\begin{remark}
Let us set $\of(t,z)=\phi(t,z) \exp(\beta A(z))$ and let $\og:\R \to \R$
be a (time-independent) test function. Then a weak formulation
of~\eqref{eq:FP'} is 
$$
\frac{d}{dt} \int_\R \of(t,z) \ \og(z) \exp(-\beta A(z)) \, dz 
= - \beta^{-1} \int_\R  \sigma^2 \
  \partial_z \of \ \partial_z \og \ \exp(-\beta A),
$$
which can be rewritten as
\begin{equation}
\label{eq:titi1}
\hspace{-2.3cm}
\frac{d}{dt} \int_{\R^n} \of(t,\xi(x)) \, \og(\xi(x)) \exp(-\beta V(x))
\, dx
= - \beta^{-1} \int_{\R^n} \nabla (\of \circ \xi) \cdot \nabla (\og \circ \xi)
\exp(-\beta V).
\end{equation}
The above weak formulation should be compared with the weak formulation of the Fokker-Planck
equation~\eqref{eq:FP_full} associated to~\eqref{eq:X}:
\begin{equation}
\label{eq:titi2}
\frac{d}{dt} \int_{\R^n} f \, g \, \exp(-\beta V) 
= - \beta^{-1} \int_{\R^n} \nabla f \cdot \nabla g \, \exp(-\beta V),
\end{equation}
where $f=\psi \exp(\beta V)$, $\psi$ is the probability distribution
function of $X_t$ satisfying~\eqref{eq:X}, and $g:\R^n \to \R$ is a
(time-independent) test function. We observe that~\eqref{eq:titi1}
is~\eqref{eq:titi2} for functions which depend on $x$ only through $\xi(x)$.
\end{remark}

We now discuss the relation between the dynamics~\eqref{eq:y} that we
propose and the dynamics~\eqref{eq:zbar}. If the function $\xi$ is such
that $|\nabla \xi|=1$, then $\sigma=1$, and in view of~\eqref{eq:A'},
\eqref{eq:F} 
and~\eqref{eq:b}, we have $b=-A'$. Hence, in this case, the
effective dynamics~\eqref{eq:y} is exactly~\eqref{eq:zbar}. 
The fact that $|\nabla \xi|=1$ is equivalent to say that $\xi$ is the
signed distance to the submanifold $\Sigma_0=\{x; \,
\xi(x)=0\}$. Examples of such reaction coordinates include $\xi(x_1,
\ldots,x_n)=x_1$, or $\xi(x)=|x|$. 

More generally, assume that $\xi$ is such that $\sigma=1$. Then, in view
of~\eqref{eq:detailed}, we have $b=-A'$, and again~\eqref{eq:y} is
exactly~\eqref{eq:zbar}. Note however that, in general, $\sigma$ is not
a constant function, and~\eqref{eq:y} differs from~\eqref{eq:zbar}. We
will confirm in Section~\ref{sec:num_D2} that~\eqref{eq:y}
and~\eqref{eq:zbar} may lead to significantly different numerical
results. 

\begin{remark}
Note that $\sigma=1$ writes
$$
\int_{\Sigma_z} \exp(-\beta V) \ |\nabla \xi| \ d\sigma_{\Sigma_z}
= 
\int_{\Sigma_z} \exp(-\beta V) \ |\nabla \xi|^{-1} \ d\sigma_{\Sigma_z}.
$$
Differentiating this equality with respect to $z$ yields (using again
Lemma~\ref{lem:chi_bar_prime}) 
\begin{eqnarray*}
\int_{\Sigma_z} \left( - \nabla V \cdot \nabla \xi  + \beta^{-1} \Delta
  \xi \right) \exp(-\beta V) \ |\nabla \xi|^{-1} \ d\sigma_{\Sigma_z}  
\\
= - \int_{\Sigma_z} \left( \frac{\nabla V \cdot \nabla \xi}{|\nabla
  \xi|^2}  - \beta^{-1} \mbox{\rm div}
\left( \frac{\nabla \xi}{|\nabla \xi|^2} \right) \right) \ \exp(-\beta
  V) \ |\nabla \xi|^{-1} \ d\sigma_{\Sigma_z},
\end{eqnarray*}
which is exactly $b=-A'$.
\end{remark}

Actually, using the fact that $\xi$ is a {\em scalar} function, it is
possible to recover the case $\sigma=1$ (for which the effective
dynamics is driven by the free energy) by two different methods. It is
not clear to us 
whether such a reformulation is also possible in the case of a
multi-dimensional reaction coordinate.

\medskip

A first method is to introduce the
following reindexation of the foliation $(\Sigma_z)_{z \in \R}$. We set
$$
h(x)=\int_0^x \sigma^{-1}(y)\,dy
$$
and we introduce the new reaction coordinate
$$
\zeta=h \circ \xi.
$$
Note that the foliation associated with $\zeta$ is exactly the same
as the one associated with $\xi$ since $h: \R \to \R$ is a one-to-one function. It
is then easy to check that the coarse-grained dynamics associated with
the reaction coordinate $\zeta$ is 
\begin{equation}
\label{eq:tata}
dy_t = - {\cal A}'(y_t) \, dt + \sqrt{2 \beta^{-1}} \, d B_t,
\end{equation}
where ${\cal A}$ is the free energy associated to $\zeta$. We hence
obtain a dynamics of the type~\eqref{eq:zbar}, with an appropriate
noise (that is, $dB_t$ in~\eqref{eq:tata} and $dW_t$ in~\eqref{eq:X} are
linked by~\eqref{eq:dB}).

\medskip

Another possibility is to keep $\xi$ as the reaction coordinate, and to
consider, instead of~\eqref{eq:X}, the dynamics
$$
dX_t = - \nabla (V - \beta^{-1}\ln (|\nabla \xi|^{-2})) \ |\nabla
\xi|^{-2}(X_t) \ dt + \sqrt{2 \beta^{-1}} \ |\nabla \xi|^{-1}(X_t) \ d W_t.
$$
The measure $\mu$ is also invariant for this dynamics. 
Then, following the same coarse-graining procedure, based on the
reaction coordinate $\xi$, one ends up with the coarse-grained dynamics
$$
dy_t = - A'(y_t) \, dt + \sqrt{2 \beta^{-1}} \ d B_t,
$$
where $A$ is the free energy associated to $\xi$. This is
exactly~\eqref{eq:zbar}, again with an appropriate noise.

\section{Error estimation in terms of time marginals}
\label{sec:preuve} 

In this section, we establish conditions on $\xi$ under which
the effective dynamics~\eqref{eq:y} is close to the dynamics of
$\xi(X_t)$, from the time marginals viewpoint ([D3] in our above
classification). 

\subsection{Error estimation}

Let $\psi^\xi(t,z)$ be the probability distribution function of $\xi(X_t)$,
where $X_t$ follows~\eqref{eq:X}, and $\phi(t,z)$ be the probability
distribution function of the solution $y_t$ to~\eqref{eq:y}. Our aim is
to bound the distance, for any time $t$, between these two
one-dimensional probability measures. 

We already introduced the total variation norm to measure distances
between measures. In the case of {\em probability
  measures}, there are two other useful quantities. The first one is the
relative entropy, which is defined by  
$$
H \left( \nu | \eta \right) = 
\int \ln\left(\frac{d\nu}{d\eta}\right) d\nu,
$$
for any two probability measures $\nu$ and $\eta$ such that $\nu$ is
absolutely continuous with respect to $\eta$. 
The relative entropy provides an upper-bound on the total variation norm
distance, by the Csisz\'ar-Kullback inequality:
\begin{equation}
\label{eq:CK}
\| \nu - \eta \|_{TV} \leq \sqrt{2 H\left( \nu | \eta \right)}.
\end{equation}
The second one is the Wasserstein distance with quadratic cost, which
is defined, for any two probability measures $\nu$ and $\eta$ with
support on a 
Riemannian manifold $\Sigma$, by 
$$
W(\nu,\eta) = \sqrt{\inf_{\pi \in \Pi(\nu,\eta)} 
\int_{\Sigma \times \Sigma} d_{\Sigma}(x,y)^2 \  d\pi(x,y)}.
$$ 
In the above expression, $d_{\Sigma}(x,y)$ denotes the geodesic distance
between $x$ and $y$ on $\Sigma$,
$$
d_{\Sigma}(x,y) = \inf \left\{ \sqrt{ \int_0^1 \left| \dot{\alpha}(t) \right|^2 \, dt};
\ \alpha \in C^1([0,1],\Sigma), \, \alpha(0) = x, \, \alpha(1) = y \right\},
$$
and $\Pi(\nu,\eta)$ denotes the set of coupling probability measures, that is
probability measures $\pi$ on $\Sigma \times \Sigma$ such that
their marginals are $\nu$ and $\eta$: for any test function $\Phi$, 
$$
\int_{\Sigma \times \Sigma} \Phi(x) \, d\pi(x,y) = \int_{\Sigma} \Phi(x)
\, d\nu(x)
\quad \text{and} \quad
\int_{\Sigma \times \Sigma} \Phi(y) \, d\pi(x,y) = \int_{\Sigma}
  \Phi(y) \, d\eta(y).
$$

In the sequel, we will need two functional inequalities, that we now
recall~\cite{ABC-00}:

\begin{defi}
A probability measure $\eta$ satisfies a logarithmic
Sobolev inequality with a constant $\rho>0$ if, for any probability
measure $\nu$, 
$$
H(\nu | \eta) \leq \frac{1}{2 \rho} I(\nu | \eta)
$$
where the Fisher information $I(\nu | \eta)$ is defined by
$$
I(\nu | \eta) = \int \left| \nabla \ln \left( \frac{d\nu}{d\eta} \right)
\right|^2 d\nu. 
$$
\end{defi}

\begin{defi}
A probability measure $\eta$ satisfies a Talagrand inequality with a
constant $\rho>0$ if, for any probability measure $\nu$,
$$
W(\nu,\eta) \leq \sqrt{ \frac{2}{\rho} H(\nu | \eta)}.
$$
\end{defi}

We will also need the following important result
(see~\cite[Theorem~1]{otto-villani-00} and~\cite{bobkov}): 
\begin{lemma}
\label{lem:lsi-t}
If $\eta$ satisfies a logarithmic Sobolev inequality with a constant
$\rho>0$, then $\eta$ satisfies a Talagrand inequality with the same
constant $\rho>0$.
\end{lemma}

Logarithmic Sobolev inequalities are very useful to prove properties
concerning the longtime behaviour of solutions to PDEs (e.g. long time
convergence of the solution of a Fokker-Planck equation to the
stationary measure of the corresponding SDE). We refer to
\cite{ABC-00,arnold-markowich-toscani-unterreiter-01,villani-03} for 
more details on this subject. 

We are now in position to present the main result of this section.

\begin{proposition}
\label{prop:D3}
Assume that $\xi$ satisfies {\bf [H1]}, and 
that the conditioned probability measures $\mu_{\Sigma_z}$,
defined by~\eqref{eq:mu_z}, satisfy a logarithmic Sobolev inequality
with a constant $\rho$ uniform in $z$: for any probability measure $\nu$
on $\Sigma_z$ which is absolutely continuous with respect to the measure
$\mu_{\Sigma_z}$, we have 
$$
  \text{{\bf [H2]}~~~}
H(\nu|\mu_{\Sigma_z}) \le \frac{1}{2 \rho} I(\nu|\mu_{\Sigma_z}) .
$$
Let us also assume that the coupling is bounded in the following sense:
$$
  \text{{\bf [H3]}~~~}
\kappa = \|\nabla_{\Sigma_z} F \|_{L^\infty} < \infty,
$$
where $F$ is the local mean force defined by~\eqref{eq:F}. 

Finally, let us assume that $|\nabla \xi|$ is close to a constant on the
manifold $\Sigma_z$ in the following sense: 
$$
\text{{\bf [H4]}~~~} 
\lambda=\left\| \frac{|\nabla \xi|^2 - \sigma^2 \circ \xi}{\sigma^2 \circ \xi}
\right\|_{L^\infty} < \infty.
$$
Assume that, at time $t=0$, the distribution of the initial conditions
of~\eqref{eq:X} and~\eqref{eq:y} are consistent one with each other:
$\psi^\xi(t=0,\cdot) = \phi(t=0,\cdot)$. 
Then we have the following estimate: for any time $t \ge 0$,
\begin{equation}
\label{eq:D3estim}
E(t)
\le \frac{M^2}{4m^2} \left( \lambda^2 + \frac{m^2 \beta^2 \kappa^2}{\rho^2}
\right) \left(
H(\psi(0,\cdot) | \mu) - H(\psi(t,\cdot) | \mu)\right),
\end{equation}
where $E(t)$ is the relative entropy of the probability distribution
function $\psi^\xi$ of $\xi(X_t)$, where $X_t$ follows~\eqref{eq:X},
with respect to the probability distribution function $\phi$ of the solution
$y_t$ to~\eqref{eq:y}:
$$
E(t)=H \left( \psi^\xi(t,\cdot) | \phi(t,\cdot) \right) = 
\int_\R \ln\left(\frac{\psi^\xi(t,z)}{\phi(t,z)}\right) \psi^\xi(t,z) \,
dz.
$$
\end{proposition}
Let us comment on these three assumptions. 
Assumption {\bf [H2]} means that $\mu_{\Sigma_z}$, which is a measure on
the manifold $\Sigma_z$, is easy to sample from. In view
of~\eqref{eq:D3estim}, the interesting case is when $\rho$ is large,
and then assumption {\bf [H2]}
implies that there is no metastability in the manifold
$\Sigma_z$. This amounts to assuming that the overdamped dynamics with respect to $\mu_{\Sigma_z}$ (which lives on $\Sigma_z$) is well-mixing. Note
finally that, in view of~\eqref{eq:mu_z_bis}, the relative entropy 
$H(\nu|\mu_{\Sigma_z})$ and the Fisher information 
$I(\nu|\mu_{\Sigma_z})$ entering assumption {\bf [H2]} read
$$
H(\nu|\mu_{\Sigma_z}) = \int_{\Sigma_z} \ln\left( f \Big/
  \frac{Z^{-1} \exp(-\beta V) }{ \exp(-\beta A(z))}\right) \ f \ |\nabla
\xi|^{-1} 
d\sigma_{\Sigma_z}
$$
and
$$
I(\nu | \mu_{\Sigma_z}) = \int_{\Sigma_z}
\left|\nabla_{\Sigma_z} \ln\left( \frac{f}{\exp(-\beta V)}
\right)\right|^2 \ f \ |\nabla \xi|^{-1} d\sigma_{\Sigma_z}, 
$$
where $f$ is the density of $\nu$ with respect to the measure $|\nabla
\xi|^{-1} \sigma_{\Sigma_z}$, {\em i.e.}
$\dis{ f = \frac{d\nu}{|\nabla \xi|^{-1} d\sigma_{\Sigma_z}}}$, and
$\nabla_{\Sigma_z}$ denotes the surface gradient: 
$$
\nabla_{\Sigma_z} = P \nabla, \quad \text{where} \quad
P(x)={\rm Id} - \frac{\nabla \xi \otimes \nabla \xi}{|\nabla \xi|^2}(x)
$$ 
is the orthogonal projector on the tangent space to $\Sigma_z$ at point $x \in \Sigma_z$.

We now turn to assumption {\bf [H3]}. Consider first the case
when $x=(x_1,x_2) \in \R^2$, and $\xi(x) = x_1$. Then $F = \nabla_{x_1} V$
and $\nabla_{\Sigma_z} F = \nabla_{x_2} F = \nabla_{x_1 x_2} V$. Requesting that
$\kappa$ is small hence amounts to requesting that $\nabla_{x_1 x_2} V$ is
small, where $x_1$ is the reaction coordinate direction whereas $x_2$ is the
direction in $\Sigma_z$. We hence ask for the coupling of these two
directions to be small. In particular, in the case when $V(x) = \frac12
x^T H x$ for some symmetric positive matrix $H \in \R^{n \times n}$ and
$\xi(x) = \xi(x_1,\ldots,x_n) = (x_1,\ldots,x_p)$ for some $p \leq n$,
we have that $\nabla_{\Sigma_z} F = 0$ if and only if the covariance
$\mbox{Cov}_\mu \left( \left( X_1,\ldots,X_p
  \right),\left(X_{p+1},\ldots,X_n \right) \right) = 0$, where $X \in \R^n$ is
distributed according to $d\mu = Z^{-1} \exp(-\beta V(x)) \, dx$. Hence
{\bf [H3]} means that the variables $\left( X_1,\ldots,X_p \right)$, which represent
the reaction coordinate directions, are decoupled from the variables 
$\left(X_{p+1},\ldots,X_n \right)$, which represent the
directions of $\Sigma_z$. 

In Section~\ref{sec:para}, we will consider an explicit example, and
compute an estimation of $\rho$ and $\kappa$ in that case, which will
help understanding the assumptions {\bf [H2]} and {\bf [H3]}.

The assumption {\bf [H4]} is technical. Observe that, if $| \nabla \xi |$
is a constant number in each manifold $\Sigma_z$, then $\lambda = 0$.  

\medskip

Before proving Proposition~\ref{prop:D3}, let us comment on the
estimate~\eqref{eq:D3estim}. Note first that this estimate 
is uniform in time. The initial conditions
for~\eqref{eq:y} and~\eqref{eq:X} are such that $\phi(t=0,\cdot) =
\psi^\xi(t=0,\cdot)$, which explains that $E(t=0) = 0$. In the longtime
limit, the estimate~\eqref{eq:D3estim} is not optimal since we know that
both $\phi$ and $\psi^\xi$ converge to $\xi 
\star \mu$ (see Lemma~\ref{lem:CG_sampling}).
This implies that $\lim_{t \to \infty} E(t) = 0$, a property that we prove in
Corollary~\ref{cor:long_time} below. 
 
\bigskip

To prove Proposition~\ref{prop:D3}, we will need the following lemma:

\begin{lemma}
\label{lem:dizdar} 
Let $\psi:\R \times \R^n \to \R$ be the probability distribution function
of $X_t$ that solves~\eqref{eq:X}. The probability distribution function of $\xi(X_t)$ is
$\dis{ \psi^\xi(t,z)=\int_{\Sigma_z} \psi(t,\cdot) |\nabla \xi|^{-1}
  d\sigma_{\Sigma_z}}$, and satisfies
\begin{eqnarray}
\hspace{-1cm}
\exp(-\beta A) \, \partial_z (\psi^\xi \exp (\beta A))
&=&
\int \frac{\nabla (\psi \exp (\beta V)) \cdot \nabla \xi}{|\nabla \xi|^2}
\ \exp(-\beta V) \ |\nabla \xi|^{-1} \ d\sigma_{\Sigma_z} 
\nonumber 
\\
& + & \beta \left(A'(z) -
  \frac{\int F \, \psi \, |\nabla \xi|^{-1} \, d\sigma_{\Sigma_z}}
{\psi^\xi}\right) \psi^\xi,
\label{eq:dizdar}
\end{eqnarray}
where $A$ is the free energy~\eqref{eq:A} and $F$ is the local mean
force~\eqref{eq:F}. 
\end{lemma}

\begin{proof}
Using~\eqref{eq:derpsi}, we compute
\begin{eqnarray*}
\exp (-\beta A) \, \partial_z (\psi^\xi \exp (\beta A))
\nonumber
\\
\quad
=
\partial_z \psi^\xi + \beta \, A' \, \psi^\xi,
\\
\quad
=
\int_{\Sigma_z} \left( 
\frac{\nabla \xi \cdot \nabla \psi}{|\nabla \xi|^2} + 
\dive \left(\frac{\nabla \xi}{|\nabla \xi|^2} \right) \psi
\right) \ |\nabla \xi|^{-1} \ d\sigma_{\Sigma_z} 
+ \beta \, A' \, \psi^\xi,
\\
\quad
=
\int_{\Sigma_z}
\frac{\nabla \xi \cdot \nabla (\psi \exp (\beta V))}{|\nabla
    \xi|^2} \ \exp (-\beta V) \ |\nabla \xi|^{-1} \ d\sigma_{\Sigma_z} 
\\
\quad \quad
+ \int_{\Sigma_z} \left( \dive\left(\frac{\nabla \xi}{|\nabla \xi|^2}
\right)  - \beta \frac{\nabla \xi \cdot \nabla V}{|\nabla
    \xi|^2}\right) \ \psi \ |\nabla \xi|^{-1} \ d\sigma_{\Sigma_z} 
+ \beta \, A' \, \psi^\xi,
\end{eqnarray*}
which yields~\eqref{eq:dizdar}.
\end{proof}

\medskip

We are now in position to prove Proposition~\ref{prop:D3}.

\medskip

\begin{proof}
We know that $\phi$ satisfies the Fokker-Planck equation~\eqref{eq:FP'},
and that $\psi^\xi$ satisfies the equation~\eqref{eq:tFP}. Thus, we have:
\begin{eqnarray*}
\frac{dE}{dt}
&=&\int \partial_t \psi^\xi \ \ln\left(\frac{\psi^\xi}{\phi}\right) 
- \int \partial_t \phi \ \frac{\psi^\xi}{\phi},
\\
&=&
\int \partial_z\left( - \tb \ \psi^\xi + \beta^{-1}
  \partial_z (\tsigma^2 \ \psi^\xi)\right) \
\ln\left(\frac{\psi^\xi}{\phi}\right)
\\
& \quad & -\beta^{-1} \int 
\partial_z\left[ \sigma^2 \, \partial_z ( \phi \exp(\beta A)) \exp(-\beta
  A) \right] 
\ \frac{\psi^\xi}{\phi},
\\
&=&- \int \left( - \tb \ \psi^\xi + \beta^{-1} \partial_z 
(\tsigma^2 \psi^\xi)\right) \partial_z \ln
\left(\frac{\psi^\xi}{\phi}\right)
\\
& \quad& + \beta^{-1} \int \sigma^2 \ \partial_z ( \phi \exp(\beta A))
\exp(-\beta A) \ \partial_z
\left(\frac{\psi^\xi}{\phi}\right).
\end{eqnarray*}
Using~\eqref{eq:dersigpsi}, we have:
\begin{eqnarray*}
\partial_z(\tsigma^2 \, \psi^\xi)
&=&
\int_{\Sigma_z} \left( \nabla \xi \cdot \nabla \psi +
\psi \Delta \xi \right) \ |\nabla \xi|^{-1} \ d\sigma_{\Sigma_z}
\\
&=&
\int_{\Sigma_z} \left( \nabla \xi \cdot \nabla (\psi \exp(\beta V))
  \exp(-\beta V) \right) \ |\nabla \xi|^{-1} \ d\sigma_{\Sigma_z}
\\
& & 
+
\int_{\Sigma_z} \left( - \beta \nabla \xi \cdot \nabla V 
 + \Delta \xi \right) \ \psi \ |\nabla \xi|^{-1} \ d\sigma_{\Sigma_z}
\\
&=&
\int_{\Sigma_z} \left( \nabla \xi \cdot \nabla (\psi \exp(\beta V)) \exp(-\beta V) \right)
\ |\nabla \xi|^{-1} \ d\sigma_{\Sigma_z} 
\\
& & + \beta \, \tb(t,z) \, \psi^\xi(t,z).
\end{eqnarray*}
Thus, it holds:
\begin{equation*}
\hspace{-2.5cm}
\begin{array}{rcl}
\dis
\frac{dE}{dt}
&=& \dis
- \beta^{-1} \int \int_{\Sigma_z} \left( 
\nabla \xi \cdot \nabla (\psi \exp(\beta V)) \exp(-\beta V) \right)
\ |\nabla \xi|^{-1} \ d\sigma_{\Sigma_z} \
\partial_z \ln \left(\frac{\psi^\xi}{\phi}\right)
\\
& \quad & \dis
+ \beta^{-1} \int \sigma^2 \ \partial_z ( \phi \exp(\beta A))
\exp(-\beta A)  \ \partial_z \left(\frac{\psi^\xi}{\phi}\right),
\\
&=& \dis
- \beta^{-1} \int \int_{\Sigma_z} 
\left( \frac{\nabla \xi \cdot \nabla (\psi \exp(\beta V))}{|\nabla \xi|^2}
\ \exp(-\beta V) \right) \, \left(|\nabla \xi|^2 - \sigma^2(z) \right) \,
\ |\nabla \xi|^{-1} \ d\sigma_{\Sigma_z} \ \partial_z \ln
\left(\frac{\psi^\xi}{\phi}\right)
\\
& \quad & \dis
-\beta^{-1} \int \sigma^2(z) \int_{\Sigma_z} \left( 
   \frac{\nabla \xi \cdot \nabla (\psi \exp(\beta V))}{|\nabla \xi|^2}
\ \exp(-\beta V) \right) \ |\nabla \xi|^{-1} \ d\sigma_{\Sigma_z} \ 
\partial_z \ln \left(\frac{\psi^\xi}{\phi}\right)
\\
& \quad & \dis
+ \beta^{-1} \int \sigma^2 \ \partial_z ( \phi \exp(\beta A))
\exp(-\beta A) \ \partial_z \left(\frac{\psi^\xi}{\phi}\right).
\end{array}
\end{equation*}
We next use~\eqref{eq:dizdar} to get:
\begin{equation*}
\hspace{-2.5cm}
\begin{array}{rcl}
\dis
\frac{dE}{dt}
&=& \dis
- \beta^{-1} \int \int_{\Sigma_z} \left( 
   \frac{\nabla \xi \cdot \nabla (\psi \exp(\beta V))}{|\nabla \xi|^2}
\ \exp(-\beta V) \right) \left(|\nabla \xi|^2 - \sigma^2(z) \right)
\ |\nabla \xi|^{-1} \ d\sigma_{\Sigma_z} \ 
\partial_z \ln \left(\frac{\psi^\xi}{\phi}\right) 
\\
& \quad & \dis
-\beta^{-1} \int \sigma^2 \left[ (\exp (-\beta A)) \, \partial_z
  (\psi^\xi \exp (\beta A)) - \beta \left(A'(z) -
   \frac{\int F \psi |\nabla \xi|^{-1}
     d\sigma_{\Sigma_z}}{\psi^\xi}\right) \psi^\xi \right] 
\ \partial_z \ln \left(\frac{\psi^\xi}{\phi}\right)
\\
& \quad & \dis
+ \beta^{-1} \int \sigma^2 \ \partial_z ( \phi \exp(\beta A))
\ \exp(-\beta A) \ \partial_z \left(\frac{\psi^\xi}{\phi}\right)
\\
&=& \dis
- \beta^{-1} \int \int_{\Sigma_z} \left( 
\frac{\nabla \xi \cdot \nabla (\psi \exp(\beta V))}{|\nabla \xi|^2}
\ \exp(-\beta V) \right) \ \left(|\nabla \xi|^2 - \sigma^2(z) \right) \
|\nabla \xi|^{-1} \ d\sigma_{\Sigma_z} 
\ \partial_z \ln \left(\frac{\psi^\xi}{\phi}\right)
\\
& \quad & \dis
+ \int \sigma^2 \left( A'(z) - 
\frac{\int F \psi |\nabla \xi|^{-1} \
    d\sigma_{\Sigma_z}}{\psi^\xi}\right) 
\ \psi^\xi \ \partial_z \ln \left(\frac{\psi^\xi}{\phi}\right)
\\
& \quad & \dis
+ \beta^{-1} \int \sigma^2 \ \exp(-\beta A) \ \partial_z
\left(\frac{\psi^\xi}{\phi}\right) \left[ \partial_z ( \phi \exp(\beta A))
-  \partial_z (\psi^\xi \exp (\beta A)) \left(\frac{\phi}{\psi^\xi}\right)
\right],
\\
&=& \dis
- \beta^{-1} \int \int_{\Sigma_z} \left( 
\frac{\nabla \xi \cdot \nabla (\psi \exp(\beta V))}{|\nabla \xi|^2}
\ \exp(-\beta V) \right) \ \left(|\nabla \xi|^2 - \sigma^2(z) \right)
\ |\nabla \xi|^{-1} \ d\sigma_{\Sigma_z} 
\ \partial_z \ln \left(\frac{\psi^\xi}{\phi}\right)
\\
& \quad & \dis
+ \int \sigma^2 \left( A'(z) -
   \frac{\int F \psi |\nabla \xi|^{-1}
     d\sigma_{\Sigma_z}}{\psi^\xi}\right) \ \psi^\xi 
\ \partial_z \ln  \left(\frac{\psi^\xi}{\phi}\right) 
- \beta^{-1} \int \sigma^2 \ \psi^\xi \ 
\left| \partial_z \ln \left(\frac{\psi^\xi}{\phi}\right) \right|^2.
\end{array}
\end{equation*}
We now use two Young inequalities, with $\eps_1>0$ and $\eps_2>0$ to be
fixed later on: 
\begin{equation}
\hspace{-2.5cm}
\begin{array}{rcl}
\dis
\frac{dE}{dt}
&\le & \dis
\frac{\beta^{-1}}{2 \eps_1} \int \left| \int_{\Sigma_z} \left( 
   \frac{\nabla \xi \cdot \nabla (\psi \exp(\beta V))}{|\nabla \xi|^2}
\ \exp(-\beta V) \right) \ \left(|\nabla \xi|^2 - \sigma^2(z) \right)
\ |\nabla \xi|^{-1} \ d\sigma_{\Sigma_z} \right|^2 
\frac{1}{\sigma^2 \psi^\xi}
\nonumber \\
& \quad & \dis
+ \frac{\beta}{2 \eps_2} \int \sigma^2 \left( A'(z) -
   \frac{\int F \psi |\nabla \xi|^{-1}
       d\sigma_{\Sigma_z}}{\psi^\xi}\right)^2 
\ \psi^\xi 
\nonumber \\
& \quad & \dis
- \beta^{-1} \left( 1 - \frac{\eps_1+\eps_2}{2} \right) \int \sigma^2 \,
\psi^\xi \, \left| \partial_z \ln \left(\frac{\psi^\xi}{\phi}\right)
\right|^2.
\label{eq:1}
\end{array}
\end{equation}
Let us first consider the second term of~\eqref{eq:1}. We write, using
{\bf [H3]}, that
\begin{eqnarray}
\left(A'(z) - \frac{\int_{\Sigma_z} F \, \psi \, |\nabla \xi|^{-1}
\ d\sigma_{\Sigma_z}}{\psi^\xi} \right)^2
&=&
\left( \int_{\Sigma_z} F \ d\mu_{\Sigma_z}- \int_{\Sigma_z} F \
  d\psi_{\Sigma_z} \right)^2
\nonumber
\\
& \leq & 
\|\nabla_{\Sigma_z} F \|_{L^\infty}^2 \
W(d\psi_{\Sigma_z},d\mu_{\Sigma_z})^2,
\label{eq:toto}
\end{eqnarray}
where $\psi_{\Sigma_z}$ is the measure $\psi(t,x)\,dx$ conditioned to 
$\xi(x)=z$:
$$
d\psi_{\Sigma_z}=\frac{\psi \, |\nabla \xi|^{-1} \, 
d\sigma_{\Sigma_z}}{\psi^\xi}.
$$
Since $\mu_{\Sigma_z}$ satisfies a logarithmic Sobolev inequality
(assumption {\bf [H2]}), it also satisfies a Talagrand inequality (see
Lemma \ref{lem:lsi-t}), hence 
$$
W(d\psi_{\Sigma_z},d\mu_{\Sigma_z})^2 
\leq 
\frac{2}{\rho} \, H(d\psi_{\Sigma_z} | d\mu_{\Sigma_z})
\leq 
\frac{1}{\rho^2} \, I(d\psi_{\Sigma_z} | d\mu_{\Sigma_z}).
$$
Gathering the above inequality with~\eqref{eq:toto}, we obtain
$$
\left(A'(z) - \frac{\int_{\Sigma_z} F \, \psi \, |\nabla \xi|^{-1} 
\, d\sigma_{\Sigma_z}}{\psi^\xi} \right)^2
\leq
\frac{\kappa^2}{\rho^2} \ I\left( d\psi_{\Sigma_z} | d\mu_{\Sigma_z} \right).
$$
Using {\bf [H1]}, we thus bound the second term of~\eqref{eq:1}:
\begin{eqnarray}
\int_\R \sigma^2 \left(A'(z) - \frac{\int_{\Sigma_z} F \, \psi \, 
|\nabla \xi|^{-1} \, d\sigma_{\Sigma_z}}{\psi^\xi} \right)^2 \psi^\xi
\nonumber
\\
\hspace{2cm}
\le
\frac{M^2\kappa^2}{\rho^2} \int_\R I\left( d\psi_{\Sigma_z} |
  d\mu_{\Sigma_z} \right) \, \psi^\xi,
\nonumber
\\
\hspace{2cm}
=
\frac{M^2\kappa^2}{\rho^2} \int_\R \int_{\Sigma_z}
\left|\nabla_{\Sigma_z} \ln\left( \frac{\psi}{\exp(-\beta V)}
\right)\right|^2 \, \psi \, |\nabla \xi|^{-1} \, d \sigma_{\Sigma_z}, 
\nonumber
\\
\hspace{2cm}
=
\frac{M^2\kappa^2}{\rho^2} \int_{\R^n}
\left|\nabla_{\Sigma_z} \ln\left( \frac{\psi}{\exp(-\beta V)}
\right)\right|^2 \psi.
\label{eq:bound1}
\end{eqnarray}
We now bound the first term of~\eqref{eq:1} using a Cauchy-Schwarz
inequality, {\bf [H4]} and {\bf [H1]}:
\begin{eqnarray}
\hspace{-1.5cm}
\int \left| \int_{\Sigma_z} \left( 
   \frac{\nabla \xi \cdot \nabla (\psi \exp(\beta V))}{|\nabla \xi|^2}
\ \exp(-\beta V) \right) \ \left(|\nabla \xi|^2 - \sigma^2(z) \right)
\ |\nabla \xi|^{-1} \ d\sigma_{\Sigma_z} \right|^2 
\frac{1}{\sigma^2 \psi^\xi}
\nonumber
\\
=
\int \left| 
\int_{\Sigma_z} \frac{\nabla \xi \cdot \nabla \ln (\psi \exp(\beta V))}{|\nabla \xi|^2}
\ \left(|\nabla \xi|^2 - \sigma^2(z) \right) \psi \ |\nabla \xi|^{-1} \ 
d\sigma_{\Sigma_z} \right|^2 \frac{1}{\sigma^2 \psi^\xi},
\nonumber
\\
\le  
\int \int_{\Sigma_z} \left| 
\frac{\nabla \xi \cdot \nabla \ln (\psi \exp(\beta V))}{|\nabla \xi|^2}
\ \left(|\nabla \xi|^2 - \sigma^2(z)\right)  \right|^2 \ \psi \
|\nabla \xi|^{-1} \ d\sigma_{\Sigma_z} \frac{1}{\sigma^2},
\nonumber
\\
\le
\lambda^2 \int \int_{\Sigma_z} \left| 
\frac{\nabla \xi \cdot \nabla \ln (\psi \exp(\beta V))}{|\nabla \xi|^2}
\right|^2 \ \psi \ |\nabla \xi|^{-1} \ d\sigma_{\Sigma_z} \ \sigma^2,
\nonumber
\\
\le
\lambda^2 M^2 \int_{\R^n}  \left| 
   \frac{\nabla \xi \cdot \nabla \ln (\psi \exp(\beta V))}{|\nabla \xi|^2}
     \right|^2  \psi.
\label{eq:bound2}
\end{eqnarray}
We infer from~\eqref{eq:1} and the bounds~\eqref{eq:bound1}
and~\eqref{eq:bound2} that
\begin{eqnarray*}
\frac{dE}{dt}
&\le& 
\frac{\beta^{-1}}{2 \eps_1} \lambda^2 M^2 \int_{\R^n}  \left| 
\frac{\nabla \xi \cdot \nabla \ln (\psi \exp(\beta V))}{|\nabla \xi|^2}
\right|^2  \psi 
\\
& \quad & +
\frac{\beta}{2 \eps_2} \frac{M^2\kappa^2}{\rho^2} \int_{\R^n}
\left| \nabla_{\Sigma_z} \ln\left( \psi \exp(\beta V) \right)\right|^2
\psi
\\
& \quad &     
- \beta^{-1} \left( 1 - \frac{\eps_1+\eps_2}{2} \right)  
\int \sigma^2 \, \psi^\xi \, \left| \partial_z 
\ln \left(\frac{\psi^\xi}{\phi}\right) \right|^2.
\end{eqnarray*}
Note that
$$
\left|\nabla \ln\left( \psi \, \exp(\beta V) \right)\right|^2
=
\left|\frac{\nabla \xi \cdot \nabla \ln (\psi \exp(\beta V))}{|\nabla \xi|}
\right|^2
+
\left|\nabla_{\Sigma_z} \ln\left( \psi \exp(\beta V) \right)\right|^2.
$$
Using the lower bound on $|\nabla \xi|$ given by {\bf [H1]}, we hence obtain
\begin{eqnarray*}
\hspace{-1.5cm}
\frac{dE}{dt}
\le
\frac{\beta^{-1}}{2 \eps_1} \frac{\lambda^2 M^2}{m^2} \int_{\R^n}  
\left| \nabla \ln (\psi \exp(\beta V)) \right|^2  \psi 
+  
\frac{\beta}{2 \eps_2} \frac{M^2\kappa^2}{\rho^2} \int_{\R^n}
\left| \nabla \ln\left( \psi \exp(\beta V) \right)\right|^2
\psi
\\
\quad
- \beta^{-1} \left( 1 - \frac{\eps_1+\eps_2}{2} \right)  
\int \sigma^2 \, \psi^\xi \, \left| \partial_z 
\ln \left(\frac{\psi^\xi}{\phi}\right) \right|^2.
\end{eqnarray*}
We now optimize on $\eps_1$ and $\eps_2$ by choosing them such that
$\eps_1+\eps_2=2$ and $\dis{ \frac{\beta^{-1}}{2 \eps_1}
\frac{\lambda^2 M^2}{m^2} = \frac{\beta}{2 \eps_2}
\frac{M^2\kappa^2}{\rho^2} }$.
This yields $\dis{ \eps_1 = \frac{2 \lambda^2 \rho^2}{\lambda^2 \rho^2 +
m^2 \beta^2 \kappa^2} }$, thus
\begin{eqnarray*}
\frac{dE}{dt}
& \le & 
\frac{\beta^{-1} M^2}{4m^2} \left( \lambda^2 + \frac{m^2 \beta^2 \kappa^2}{\rho^2}
\right) \int_{\R^n} \left| \nabla \ln (\psi \exp(\beta V)) \right|^2 \psi,
\\
&=&  
\frac{\beta^{-1} M^2}{4m^2} \left( \lambda^2 + \frac{m^2 \beta^2 \kappa^2}{\rho^2}
\right) I(\psi|\mu),
\\
&=&
-\frac{M^2}{4m^2} \left( \lambda^2 + \frac{m^2 \beta^2 \kappa^2}{\rho^2}
\right) \, \frac{d}{dt} H(\psi|\mu).
\end{eqnarray*}
We next integrate this equation between $0$ and $t$ and use the fact
that $E(0)=0$ to obtain~\eqref{eq:D3estim}.
\end{proof}

\medskip

We now prove a corollary of Proposition~\ref{prop:D3}, which strengthens
its long-time limit behaviour.

\begin{coro}
\label{cor:long_time}
In addition to the assumptions of Proposition~\ref{prop:D3}, assume that
$$
\text{{\bf [H5]}~The measure $\xi \star \mu$ satisfies a logarithmic Sobolev inequality
with a constant $r$.}$$

Consider again 
the probability distribution function $\psi^\xi$ of $\xi(X_t)$, where
$X_t$ follows~\eqref{eq:X}, and the probability distribution function
$\phi$ of the solution $y_t$ to~\eqref{eq:y}. They satisfy: 
\begin{equation}
\label{tutu0}
\forall t \ge 0, \quad
\| \psi^\xi(t,\cdot) - \phi(t,\cdot) \|_{TV} \leq \min \left( C_1(t), 
2 C_2 \, \exp(- R \, \beta^{-1} \, t)
\right),
\end{equation}
for some positive constant $R$, where 
\begin{eqnarray}
\label{eq:C_1}
C_1(t) &=& 
\sqrt{\frac{M^2}{2m^2} \left( \lambda^2 + \frac{m^2 \beta^2 \kappa^2}{\rho^2}
\right) \left( H(\psi(0,\cdot) | \mu) - H(\psi(t,\cdot) | \mu) \right)},
\\
\label{eq:C_2}
C_2 &=& \max \left( \sqrt{2 H\left( \phi(0,\cdot) | \mu^\xi \right)}, 
\sqrt{2 H\left( \psi(0,\cdot) | \mu \right)} \right),
\end{eqnarray}
with $d\mu^\xi = \exp(-\beta A(z)) \, dz$.
\end{coro}
As a consequence of this corollary, we see that 
$\lim_{t \to \infty} \| \psi^\xi(t,\cdot) - \phi(t,\cdot) \|_{TV} = 0$. 

\medskip

\begin{proof}
We infer from the Csisz\'ar-Kullback inequality and from the 
bound~\eqref{eq:D3estim} that
\begin{equation}
\label{tutu4}
\| \psi^\xi - \phi \|_{TV} \leq \sqrt{2 H\left( \psi^\xi | \phi \right)}
\leq C_1(t),
\end{equation}
where $C_1(t)$ is given by~\eqref{eq:C_1}. We also have
\begin{equation}
\label{tutu1}
\| \psi^\xi - \phi \|_{TV} 
\leq 
\left\| \psi^\xi - \mu^\xi \right\|_{TV} + \left\| \phi - \mu^\xi \right
\|_{TV},
\end{equation}
where $d\mu^\xi = \exp(-\beta A(z)) \, dz$ is the equilibrium measure
$\xi \star \mu$. Let us first upper-bound 
$\dis E_{\rm CG}(t) = H\left( \phi | \mu^\xi \right)=
\int_\R \ln\left(\frac{\phi}{\exp(-\beta A)}\right)
\phi$. Using~\eqref{eq:FP'}, we compute
\begin{eqnarray*}
\frac{dE_{\rm CG}}{dt}
&=&
\int_\R \partial_t \phi \ \ln \left(\frac{\phi}{\exp(-\beta A)}\right)
\\
&=&
\beta^{-1} \int_\R \partial_z \left[ \sigma^2 \
  \partial_z (\phi \exp(\beta A) ) \ \exp(-\beta A) \right]
\ \ln \left(\frac{\phi}{\exp(-\beta A)}\right)
\\
&=&
- \beta^{-1} \int_\R \left[ \sigma^2 \
  \partial_z (\phi \exp(\beta A) ) \ \exp(-\beta A) \right]
\ \partial_z \left[ \ln \left(\frac{\phi}{\exp(-\beta A)}\right) \right]
\\
&\leq &
- m^2 \beta^{-1} \int_\R \phi \ \left| \partial_z \left[ 
\ln \left(\frac{\phi}{\exp(-\beta A)}\right) 
\right] \right|^2
\\
&= &
- m^2 \beta^{-1} I(\phi | \mu^\xi),
\end{eqnarray*}
where we have used that $\sigma^2 \geq m^2$, which is a consequence of
{\bf [H1]} and~\eqref{eq:sigma}. Since $\mu^\xi$ satisfies a logarithmic
Sobolev inequality with constant $r$, we infer from the above bound that
$\dis 
\frac{dE_{\rm CG}}{dt} \leq - 2 \, r \, m^2 \, \beta^{-1} E_{\rm CG}.
$
Using a Gronwall lemma, we obtain
$$
H\left( \phi | \mu^\xi \right) = 
E_{\rm CG}(t) \leq E_{\rm CG}(t=0) \, \exp(- 2 \, r \, m^2 \, \beta^{-1} \, t)
= H\left( \phi(0,\cdot) | \mu^\xi \right) \, 
\exp(- 2 \, r \, m^2 \, \beta^{-1} \, t),
$$
and the Csisz\'ar-Kullback inequality then yields
\begin{equation}
\label{tutu2}
\left\| \phi - \mu^\xi \right\|_{TV}
\leq 
\sqrt{2 H\left( \phi | \mu^\xi \right)}
\leq
\sqrt{2 H\left( \phi(0,\cdot) | \mu^\xi \right)} \, 
\exp(- r \, m^2 \, \beta^{-1} \, t).
\end{equation}
We now turn to the term $\left\| \psi^\xi - \mu^\xi \right\|_{TV}$. 
For any function $\chi : \R^n \to \R$, define 
$\dis \chi^\xi(z)=\int_{\Sigma_z} \chi \, |\nabla \xi|^{-1} \,
d\sigma_{\Sigma_z}$, and observe that
$$
\int_{\R^n} | \chi(x) | \, dx
=
\int_\R \int_{\Sigma_z} \frac{| \chi |}{|\nabla \xi|} 
\, d\sigma_{\Sigma_z} \, dz
\geq
\int_\R \left|
\int_{\Sigma_z} \frac{\chi}{|\nabla \xi|} 
\, d\sigma_{\Sigma_z} \right| \, dz
=
\int_\R \left| \chi^\xi \right| \, dz
$$
which also reads 
$\left\| \chi \right\|_{TV} \geq \left\| \chi^\xi \right\|_{TV}$. We
apply this inequality with $\chi = \psi - \mu$: 
$$
\left\| \psi^\xi - \mu^\xi \right\|_{TV} \leq 
\left\| \psi - \mu \right\|_{TV} \leq 
\sqrt{2 H\left( \psi | \mu \right)}.
$$
Since $\mu^\xi$ and the conditional measures $\mu_{\Sigma_z}$ satisfy a
logarithmic Sobolev inequality (see {\bf [H5]} and {\bf [H2]}), and under
assumption {\bf [H3]}, we obtain that the measure $\mu$ also satisfies a
logarithmic Sobolev inequality with some constant $R>0$
(see~\cite{tony_jfa}). Hence, by a computation similar to the one on
$E_{\rm CG}$, we obtain 
$$
H\left( \psi | \mu \right) \leq
H\left( \psi(0,\cdot) | \mu \right) \, 
\exp(- 2 \, R \, \beta^{-1} \, t),
$$
hence
\begin{equation}
\label{tutu3}
\left\| \psi^\xi - \mu^\xi \right\|_{TV} \leq 
\sqrt{2 H\left( \psi(0,\cdot) | \mu \right)} \, 
\exp(- R \, \beta^{-1} \, t).
\end{equation}
Gathering~\eqref{tutu1}, \eqref{tutu2} and \eqref{tutu3}, we
obtain
$$
\| \psi^\xi - \phi \|_{TV} 
\leq 
C_2 \, \exp(- r \, m^2 \, \beta^{-1} \, t)
+
C_2 \, \exp(- R \, \beta^{-1} \, t),
$$
where $C_2$ is defined by~\eqref{eq:C_2}. The proof of~\cite[Theorem
1.2]{tony_jfa} shows that $0 < R \leq r m^2$. The above bound then yields
$\| \psi^\xi - \phi \|_{TV} \leq 2 C_2 \, \exp(- R \, \beta^{-1} \, t)$,
which, gathered with~\eqref{tutu4}, yields~\eqref{tutu0}.
\end{proof}

\subsection{Estimation of the upper-bound constants
  of~\eqref{eq:D3estim} in a particular case} 
\label{sec:para}

In this section, we give a very formal argument to estimate the
constants $\rho$ and $\kappa$ entering the bound~\eqref{eq:D3estim}, in
a specific case. 
Potential energies in molecular dynamics are often the sum of several
terms, with different stiffness. For instance, the potential energy of
an alkane chain, in the United Atom model~\cite{RB78}, reads
$$
V(X) = \sum_i V_2(d_{i,i+1}) + \sum_i V_3(\theta_i) + \sum_i V_4(\phi_i)
+ V_{\rm non-bonded}(X),
$$
where $d_{i,i+1}$ is the distance between atoms $i$ and $i+1$,
$\theta_i$ is the bond angle made by atoms $i-1$, $i$ and $i+1$, whereas
$\phi_i$ is the dihedral angle defined by the atoms $i+j$,
$j=-1,\ldots,2$. In general, $V_2$ is a much stiffer potential than
$V_3$, which is itself much stiffer than $V_4$.

A simple toy-model for such potential energies is 
\begin{equation}
\label{eq:ex_V}
V_\eps(X) = V_0(X) + \frac{1}{\eps} q^2(X),
\end{equation}
where $V_0$ and $q$ are two scalar-valued functions that do not depend
on the small parameter $\eps$ (see Equation~\eqref{eq:dw} and
Figure~\ref{fig:dw} below for a precise example of type~\eqref{eq:ex_V}).
For simplicity, we assume here that the reaction coordinate $\xi$ does
not depend on $\eps$, and that it is constant on
the manifolds $\Sigma_z$ (in assumption {\bf [H4]}, $\lambda =
0$). Since the relative entropy is always non-negative, the
estimate~\eqref{eq:D3estim} reads 
$$
E(t)
\le \frac{M^2}{4} \frac{\beta^2 \kappa_\eps^2}{\rho_\eps^2} \
H(\psi(0,\cdot) | \mu_\eps).
$$
We also assume that the initial condition of~\eqref{eq:X} is well
adapted to the Boltzmann measure~$\mu_\eps$, in the sense that
$H(\psi(0,\cdot) | \mu_\eps)$ is upper-bounded by a constant independent
of $\eps$. Thus the above bound reads
$$
E(t)
\le C \frac{\kappa_\eps^2}{\rho_\eps^2}
$$
for some constant $C$ independent of $\eps$. 
Our aim is to roughly estimate the coefficients $\rho_\eps$
and $\kappa_\eps$ in terms of $\eps$.   

\medskip

Since $\eps$ is small, the Boltzmann measure~\eqref{eq:mu} concentrates
on the manifold where $q=0$, and locally looks like a Gaussian measure of
variance $\eps$ around that manifold. The same holds for $\mu_{\Sigma_z}$, that
is assumed to satisfy a logarithmic Sobolev inequality (assumption {\bf [H2]}). Hence, we
typically have
$\rho_\eps = O\left( 1/\eps \right)$. 

We now compute the local mean force, defined by~\eqref{eq:F}:
\begin{eqnarray*}
F &=& \frac{\nabla V_\eps \cdot \nabla \xi}{|\nabla \xi|^2}  - \beta^{-1} \dive
\left( \frac{\nabla \xi}{|\nabla \xi|^2} \right)
\\
&=&
\frac{2}{\eps} \ q \ \frac{\nabla q \cdot \nabla \xi}{|\nabla \xi|^2} 
\ + \  
\frac{\nabla V_0 \cdot \nabla \xi}{|\nabla \xi|^2}  - \beta^{-1} \dive
\left( \frac{\nabla \xi}{|\nabla \xi|^2} \right).
\end{eqnarray*}
Recall that $\xi$ does not depend on $\eps$. If $\nabla q
\cdot \nabla \xi \neq 0$, then $F$ is of order $O\left( 1/\eps \right)$,
and so is $\kappa_\eps$. On the contrary, if $\nabla q \cdot \nabla \xi = 0$,
then $F$ is of order $O(1)$ with respect to $\eps$, and so is
$\kappa_\eps$.  

\medskip

Let us summarize our discussion. In the case when
$\nabla q \cdot \nabla \xi = 0$, it turns out that
$\rho_\eps$ is of order $1/\eps$, while $\kappa_\eps$ is of order $1$, and 
the estimate~\eqref{eq:D3estim} reads
$$
E(t) \le C \eps^2 
$$
for some constant $C$ that does not depend on $\eps$. Hence, as $\eps$
decreases to 0, the effective dynamics~\eqref{eq:y} becomes more
accurate, in the sense of [D3].
In the case when $\nabla q \cdot \nabla \xi \neq 0$,
both $\rho_\eps$ and $\kappa_\eps$ are of order $1/\eps$, and the
estimate~\eqref{eq:D3estim} reads 
$$
E(t) \le C 
$$
for some constant $C$ that does not depend on $\eps$. So the effective
dynamics~\eqref{eq:y} is not particularly accurate.
In the next section, we numerically confirm that the criterion 
\begin{equation}
\label{eq:criterion}
\nabla \xi \cdot \nabla q = 0
\end{equation} 
has indeed a significant impact on the accuracy of the effective
dynamics.

\section{Numerical results: residence time estimation}
\label{sec:num_D2}

Our aim here is twofold. First, we want to check the accuracy
of~\eqref{eq:y} in a sense related to [D2], on a simple system, and also compare this effective dynamics with the
coarse-grained dynamics~\eqref{eq:zbar} based on the free
energy. Second, we wish to assess the relevance of the
criterion~\eqref{eq:criterion}. It
seems to be an important condition for estimates in the sense of [D3] to
be meaningful. Is it also a necessary and sufficient condition in order
to obtain accurate dynamical properties~?

In the following numerical tests, we focus on the residence times. 
We have indeed already underlined that the characteristic behaviour of the
dynamics~\eqref{eq:X} is to sample a given well of the potential energy,
then suddenly hopes to another basin, and start over. Consequently, an
important quantity is the residence time that the system spends in the
well, before 
going to another one. In this section, we describe a numerical example
where we have studied such quantities, which contain dynamical
information, and are related to the estimator [D2].

\begin{figure}
\centerline{
\input{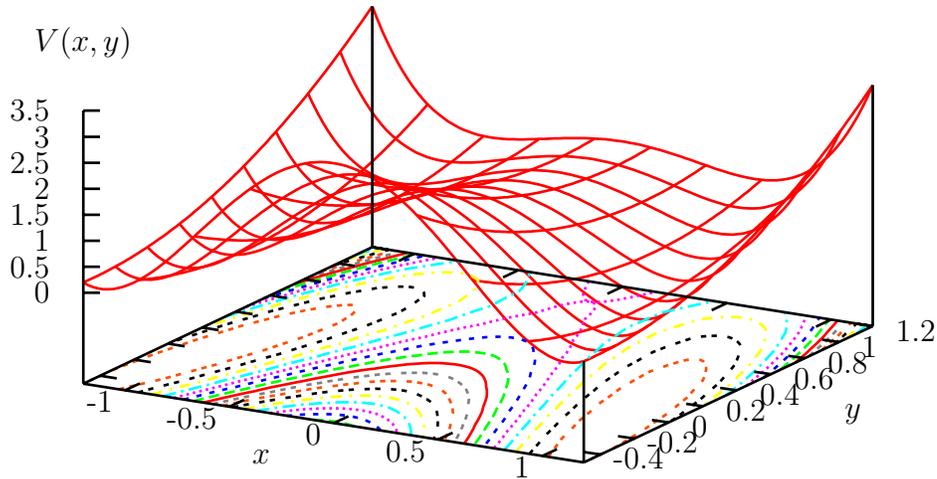}
}
\caption{Plot of the double-well potential~\eqref{eq:dw}. For clarity of
  the picture, we set $\eps = 1$.}
\label{fig:dw}
\end{figure}

Consider the two-dimensional potential energy
\begin{equation}
\label{eq:dw}
V_\eps(x,y) = (x^2 - 1)^2 + \frac{1}{\eps} (x^2 + y - 1)^2
\end{equation}
which is of the form~\eqref{eq:ex_V}, with $V_0(x,y) = (x^2 -
1)^2$ and $q(x,y) = x^2 + y - 1$. For any $\eps > 0$, the potential
$V_\eps$ has two local 
minima, at $(x,y) = (\pm 1,0)$, and one saddle point, at $(x,y) =
(0,1)$ (see Figure~\ref{fig:dw}). There are thus two basins,
namely $\left\{ (x,y) \in \R^2; \ x < 
  0 \right\}$ and $\left\{ (x,y) \in \R^2; \ x > 0 \right\}$. Since $V_\eps$
is an even function of $x$, the residence times in each
well are equal to each other. Our aim is to compare the residence time 
computed when the full description of the system is used (that is, we
simulate the dynamics~\eqref{eq:X}) with the residence time computed
from a coarse-grained description, according to~\eqref{eq:y}
or~\eqref{eq:zbar}, for two different reaction coordinates.  

In the case at hand, a natural reaction coordinate is $\xi_1(x,y) = x$,
since the value of $\xi_1$ already gives the information that the system
is in the right or the left well. In that case, $| \nabla \xi_1 | =1$,
hence the effective dynamics~\eqref{eq:y} is the same as the
dynamics~\eqref{eq:zbar}, that is the dynamics driven by the free
energy $A_1$ associated to $\xi_1$. This free energy reads
\begin{equation}
\label{eq:A1}
A_1(z) = (z^2 - 1)^2 + C(\beta)
\end{equation}
for some constant $C(\beta)$ ensuring that 
$\dis{ \int_\R \exp(-\beta A_1(z)) \, dz = 1}$. 

Note that $\nabla \xi_1 \cdot \nabla q \neq 0$. In view of the
discussion of the previous section, we do not expect the effective
dynamics based on $\xi_1$ to be very accurate. 

Consider now the function $\xi_2(x,y) = x \exp(-2y)$, which satisfies 
$\nabla \xi_2 \cdot \nabla q = 0$. We expect the effective
dynamics~\eqref{eq:y}, based on $\xi_2$, to be accurate, at least in the
sense of the estimator [D3] (time marginals). Here, we want to check its
accuracy in terms of residence times (and hence in a way related to
estimator~[D2]). Note that, for this reaction
coordinate, $| \nabla \xi_2 |$ is not a constant function, hence the
effective dynamics~\eqref{eq:y} differs from the
dynamics~\eqref{eq:zbar} for $A \equiv A_2$, the free energy associated
to $\xi_2$. 

\medskip

\begin{figure}
\centerline{
\input{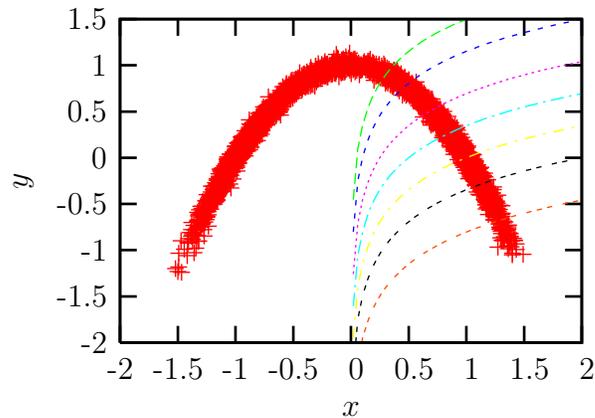}
}
\caption{Crosses: plot of the trajectory $X_t = (x_t,y_t)$ solution
  to~\eqref{eq:X}, for the parameters $\eps = 0.01$ and $\beta
  = 3$. Dashed lines: level sets of $\xi_2$.}  
\label{fig:traj_full}
\end{figure}

We work with the parameters $\eps = 0.01$ and $\beta = 3$. On
Figure~\ref{fig:traj_full}, we plot the trajectory solution
to~\eqref{eq:X}, as well as the level sets of $\xi_2$. We can see that
the trajectory remains close to the line $\{ (x,y) ; \ q(x,y) = x^2 + y
- 1 = 0 \}$
(since $\eps$ is small), and that the level sets of $\xi_2$ are parallel
to $\nabla q$, which implies that $\nabla \xi_2$ is indeed perpendicular
to~$\nabla q$. 

\begin{figure}
\centerline{
\input{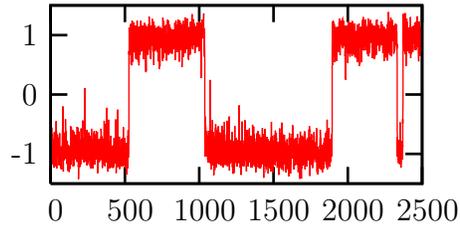}
}
\caption{Time evolution $t \mapsto x_t$, for $X_t = (x_t,y_t)$
  solution to~\eqref{eq:X}, for the parameters $\eps = 0.01$ and $\beta
  = 3$. We clearly see metastability.}  
\label{fig:traj_full_x}
\end{figure}

With the choice we made for $\beta$ and $\eps$, the system is 
metastable. On
Figure~\ref{fig:traj_full_x}, we plot $x_t$ as a function of time, where
$X_t = (x_t,y_t)$ satisfies~\eqref{eq:X}. We clearly see that $x_t$
remains close to -1 (that is, the system is in the left well) for a long
time before hoping to the right well.

\begin{figure}
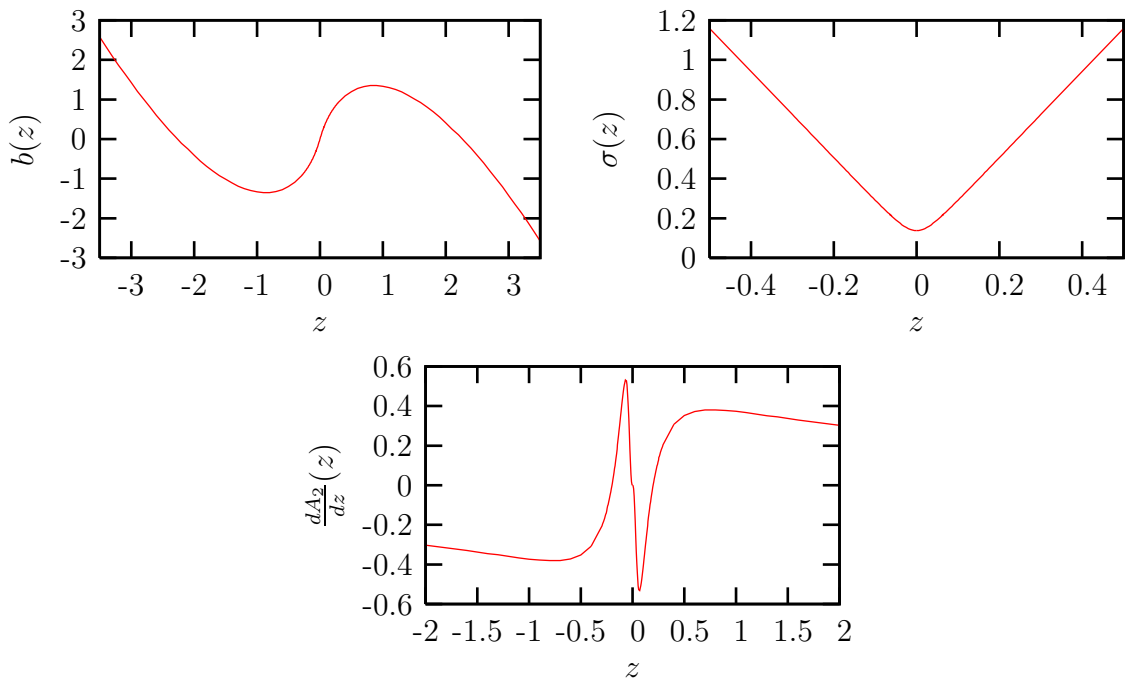

\centerline{
\input{figure4a.tex} \input{figure4b.tex}
}
\hspace{4cm}
\input{figure4c.tex}
\caption{Plot of the functions $b$, $\sigma$ and $A'_2$, for the
  reaction coordinate $\xi_2$. Note that $b$ and $A'_2$ are odd
  functions, whereas $\sigma$ is an even function. Note the large
  variations of $A'_2$ in the neighbourhood of $z=0$.}
\label{fig:A2}
\end{figure}

The functions $b$ and $\sigma$, as well as the derivative of the free
energy $A_2$  
(respectively defined by~\eqref{eq:b}, \eqref{eq:sigma}
and~\eqref{eq:A'}) are plotted on Figure~\ref{fig:A2}, in the case of the
reaction coordinate $\xi_2$. 

\begin{remark}
For all the numerical tests reported in this article, 
the complete dynamics~\eqref{eq:X} has been integrated with
the Euler-Maruyama scheme
$$
X_{j+1} = X_j -  \Delta t \, \nabla V(X_j) + \sqrt{2 \, \Delta t \,
  \beta^{-1}} \ G_j,
$$
where, for any $j$, $G_j$ is a two-dimensional vector, whose coordinates are
independent and identically distributed (i.i.d.) random variables,
distributed according to a normal Gaussian law. 

For the reaction coordinate $\xi_1$, the effective dynamics
is~\eqref{eq:zbar}, that we have numerically simulated with the same
algorithm as above. We have used the analytical expression~\eqref{eq:A1}
of the free energy $A_1$. 

For the reaction coordinate $\xi_2$, the free energy derivative $A'_2$ and the
functions $b$ and $\sigma$ have been computed using the algorithm
proposed in~\cite{c_l_eve}. We have chosen to work in the interval
$\xi_2 \in [-200;200]$, and computed $A'_2$, $b$ and $\sigma$ on a grid of size
$\Delta z = 0.1$ (except in the interval $[-0.3;0.3]$, where we used a
finer grid of size $\Delta z = 5. \, 10^{-3}$, since the variations of $A'_2$, $b$ and
$\sigma$ are larger in the neighbourhood of 0). Values of the
functions for $z$ in-between points of that grid have been obtained by linear
interpolation (see Figure~\ref{fig:A2}). We have again used the Euler-Maruyama
scheme to numerically integrate the dynamics~\eqref{eq:y}. 

All dynamics have been integrated with the time step $\Delta t = 10^{-4}$.
\end{remark}

For the reaction coordinate $\xi_i$, $i=1,2$, the left and the right
wells are defined as the sets $\left\{(x,y) \in \R^2; \ \xi_i(x,y) \leq -
  \xi_i^{\rm th} \right\}$ and $\left\{(x,y) \in \R^2; \ \xi_i(x,y) \geq 
  \xi_i^{\rm th} \right\}$, respectively. We have chosen the threshold
values $\xi_1^{\rm
  th} > 0$ and $\xi_2^{\rm th} > 0$ such that wells are more or less the same
for both reaction coordinates. 
To compute the residence time, we proceeded as follow, for both reaction
coordinates $\xi_1$ and $\xi_2$: 
\begin{enumerate}
\item we first generated $15 \, 000$ configurations 
$\{ (x_i,y_i) \in \rens^2\}_{1 \leq i \leq 15 \, 000}$, distributed
according to the measure $\mu$, and 
  such that $\xi(x_i,y_i)$ belongs to the right well, that is
  $\xi(x_i,y_i) > \xi^{\rm th}$.
\item we next ran the dynamics~\eqref{eq:X} from the initial condition
  $(x_i,y_i)$, and monitor the first time $\tau_i$ at which the
  system reaches a point $(x(\tau_i),y(\tau_i))$ in the left well:
$\tau_i = \inf \left\{ t; \ \xi(x(t),y(t)) < - \xi^{\rm th} \right\}$.
\item from these $(\tau_i)_{1 \leq i \leq 15 \, 000}$, we computed an average
  residence time and a confidence interval. These figures are
  the reference figures. 
\item we next consider the initial conditions 
$\left\{ \xi(x_i,y_i) \in \R\right\}_{1 \leq i \leq 15 \, 000}$
for the effective dynamics. By construction, these
configurations are distributed according to the equilibrium measure of
$\xi \star \mu$, that is $\exp(-\beta A(z)) \, dz$, and are in the right well. 
\item from these initial conditions, we run the dynamics~\eqref{eq:y}
  or~\eqref{eq:zbar}, until the left well is reached ($y(t) \leq -
  \xi^{\rm th}$ when working with~\eqref{eq:y}, $\zbar(t) \leq -
  \xi^{\rm th}$ when working with~\eqref{eq:zbar}), and compute, as 
  for the complete description, a residence time and its confidence
  interval.
\end {enumerate}

The results we found for the residence time are gathered in
Table~\ref{tab:residence}. We see that, when we work with $\xi_2$ (which
satisfies the condition $\nabla \xi_2 \cdot \nabla q = 0$) {\em and}
with the
effective dynamics~\eqref{eq:y}, we can reproduce the reference
residence time (32.5 $\pm$ 0.5)
within an excellent accuracy. If we still use the reaction
coordinate $\xi_2$, but consider as the coarse-grained dynamics the
dynamics~\eqref{eq:zbar} driven by the free energy $A_2$, then we obtain
results that are inconsistent with the reference results given by the complete
description of the system. 

Note also that the results obtained with choosing $\xi_1$ as reaction
coordinate, which is such that $\nabla \xi_1 \cdot \nabla q \neq 0$, are
inconsistent with the reference results (in that case,
the effective dynamics~\eqref{eq:y} is the same
as~\eqref{eq:zbar}). Actually, the coarse-grained dynamics does not
depend on $\eps$ (since the free energy $A_1$ does not depend on
$\eps$), whereas the complete description does depend on $\eps$. 

\begin{table}
\centerline{
\begin{tabular}{| c | c | c | c |c |}
  \hline
Reac. Coord. & $\xi^{\rm th}$ & Ref. residence time & 
Reduced dyn. type & CG residence time 
\\
\hline
$\xi_2(x,y)$ & 0.13 & 32.5 $\pm$ 0.5 & Eff. dyn.~\eqref{eq:y} &
32.7 $\pm$ 0.5
\\
\hline
$\xi_2(x,y)$ & 0.13 & 32.5 $\pm$ 0.5 & Dyn.~\eqref{eq:zbar} &
6.4 $\pm$ 0.3
\\
\hline \hline
$\xi_1(x,y)$ & 0.5 & 31.6 $\pm$ 0.5 & Dyn. \eqref{eq:y} = \eqref{eq:zbar}
& 24.4 $\pm$ 0.4
\\
\hline
\end{tabular}
}
\caption{Residence times obtained from the complete description (third
  column) and from the reduced description (last column), for both
  reaction coordinates (and both dynamics~\eqref{eq:y}
  and~\eqref{eq:zbar} when applicable). The threshold values
  ($\xi_1^{\rm th} = 0.5$ and $\xi_2^{\rm th} = 0.13$) have been adjusted
so that the reference residence times for both reaction coordinates
(31.6 and 32.5, respectively) are almost equal.} 
\label{tab:residence}
\end{table}

\section{Pathwise convergence}
\label{sec:path}

In this section, we prove pathwise convergence
results between $\xi(X_t)$, where $X_t$ solves~\eqref{eq:X}, and $y_t$
which solves~\eqref{eq:y}, for some potential energies of the
type~\eqref{eq:ex_V}. On these specific examples, we obtain stronger
convergence results than in the previous sections (namely, convergence
in the sense of [D1] rather than in the sense of [D3] or [D2], as in
Sections~\ref{sec:preuve} and~\ref{sec:num_D2}).

Consider the dynamics~\eqref{eq:X}, with the potential energy $V_\eps$
defined by~\eqref{eq:ex_V}. It reads
\begin{equation}
\label{eq:X_eps}
dX^\eps_t = - \nabla V_0(X^\eps_t) \, dt 
- \frac{1}{\eps} \nabla (q^2)(X^\eps_t) \, dt +
\sqrt{2 \beta^{-1}} \, d W_t,
\quad 
X^\eps_{t=0} = X_0.
\end{equation}
Note that the initial condition is supposed to not depend on $\eps$.
The limit of $X_t^\eps$ when $\eps \to 0$ has been identified in
\cite{c_l_eve}: it is a process $\overline{X}_t$ solution of a SDE that we
write below (see equation~\eqref{eq:Xbar}), and that is such that
$q(\overline{X}_t) = 0$ for any $t$. 

Assume now that $X_t^\eps \in \R^2$: then $\overline{X}_t$ belongs to the
one-dimensional manifold 
\begin{equation}
\label{eq:M}
{\cal M} = \left\{ X \in \R^2; \ q(X) = 0 \right\}.
\end{equation} 
Assume also that the reaction coordinate $\xi$ is such that
its restriction $\xi_{| {\cal M}}$ on ${\cal M}$ is a one-to-one map
from ${\cal M}$ to some subset of $\R$ (that is, $\xi$ parameterizes ${\cal
  M}$). In that case, it is easy to
build a reduced dynamics from~\eqref{eq:X_eps}, in the limit $\eps \to
0$: one first lets $\eps$ go to zero, writes the dynamics of
$\overline{X}_t$, and then makes a one-to-one change of variable to write
the dynamics in term of $\xi\left(\overline{X}_t\right)$. Our aim is to
write conditions 
under which the so-obtained dynamics corresponds to~\eqref{eq:y}, which
amounts to say that the diagram~\eqref{eq:comm} is a commutative
diagram. 

\begin{equation}
\label{eq:comm}
\hspace{-2.5cm}
\begin{array}{ccc}
\text{2D dynamics~\eqref{eq:X_eps} on } X^\eps_t
&
\dis{ \to_{\eps \to 0}^{\left[ \begin{array}{c} \text{\scriptsize pathwise}
\\ \text{\scriptsize convergence} \end{array} \right] } }
&
\text{1D limit dynamics~\eqref{eq:Xbar} on } \overline X_t
\\ 
\downarrow
&
&
\downarrow
\\
\text{\em It\^o computation}
&
&
\downarrow
\\ 
\downarrow
&
&
\downarrow
\\
\text{Nonclosed dynamics on } \xi(X^\eps_t)
&
&
\text{\em One-to-one change of variable:} 
\\
\downarrow
&
&
z_t = \xi(\overline X_t)
\\
\text{\em Conditional expectations}
&
&
\downarrow
\\
\downarrow
&
&
\downarrow
\\
\text{Dynamics~\eqref{eq:y}
on } y^\eps_t \approx \xi(X^\eps_t):
&
\to_{\eps \to 0}
&
\text{Dynamics~\eqref{eq:path_cv} on } z_t
\\
dy^\eps_t = b_\eps(y^\eps_t) \, dt + \sigma_\eps(y^\eps_t) dB_t
\end{array}
\end{equation}

\subsection{Limit of~\eqref{eq:X_eps} in a pathwise sense}

We now proceed in details. For any $X \in {\cal M}$, let 
$$
P(X) = \mbox{Id} - \frac{\nabla q \otimes \nabla q}{| \nabla q |^2}(X)
$$
be the projector on ${\cal T}_X {\cal M}$, the tangent space to ${\cal
  M}$ at $X$. Let us define
$$
n = \frac{\nabla q}{| \nabla q |} \quad \mbox{and} \quad \kappa = \dive n.
$$
Let us now introduce the process $\overline{X}_t$ solution to the
equation
\begin{equation}
\label{eq:Xbar}
\hspace{-2.3cm}
d\overline{X}_t = - P \left( \overline{X}_t \right) \nabla 
\left( V_0 + \beta^{-1} \ln | \nabla q | \right) \left( \overline{X}_t
\right) \, dt
- \beta^{-1} \kappa \, n \, dt + \sqrt{2 \beta^{-1}} \, 
P \left( \overline{X}_t \right) \, dW_t,
\end{equation}
with the same initial condition $\overline{X}_{t=0} = X_0$
as~\eqref{eq:X_eps}. 
Let us assume that this initial condition
satisfies $q(X_0)=0$, and let us fix a time interval $[0,T]$. Then
(see~\cite{c_l_eve}), under some regularity assumptions on $q$ and $V_0$,
there exists a constant $C$ that does not depend on $\eps$ such that 
\begin{equation}
\label{eq:estim_tony}
\sup_{t \in [0,T]} \E
\left| X_t^\eps - \overline{X}_t \right|^2 \leq C \eps.
\end{equation}
Note also that $q\left( \overline{X}_t \right) = 0$ for any time $t$.  

Assume now that there exists a one-to-one map
\begin{equation}
\label{eq:def_chi}
\chi: X \in \rens^2 \mapsto (\xi(X),q(X)),
\end{equation}
which implies that the manifold ${\cal M}$ defined by~\eqref{eq:M} can
be parameterized by $\xi$. Then, equation~\eqref{eq:Xbar} is equivalent
to the dynamics 
$$ 
d \left( \xi \left( \overline{X}_t \right) \right)
= 
\nabla \xi\left( \overline{X}_t \right) \cdot d\overline{X}_t + 
\beta^{-1} \, 
P\left( \overline{X}_t \right) : \nabla^2 \xi\left( \overline{X}_t
\right) \ dt.
$$
After some tedious
but not difficult computations, we see that the above dynamics can be
written 
\begin{equation}
d \left( \xi \left( \overline{X}_t \right) \right)
=
d_1\left( \overline{X}_t \right) \, dt + 
d_2\left( \overline{X}_t \right) \, dt + 
\sqrt{2 \beta^{-1}} \, | \nabla \xi | \left( \overline{X}_t \right) \, dB_t,
\label{eq:titi}
\end{equation}
with again $\dis{ dB_t = \frac{\nabla \xi}{|\nabla \xi|} \left( \overline{X}_t
  \right) \cdot dW_t}$, and 
\begin{eqnarray}
\nonumber
d_1 &=& - \nabla \xi \cdot \nabla V_0 + \beta^{-1} \, \Delta \xi,
\\
\label{eq:d2}
d_2 &=& 
- \frac{1}{\beta} \frac{\nabla q \cdot \nabla u}{| \nabla q |^2}
+ u \, \frac{\nabla q \cdot \nabla V_0}{| \nabla q |^2}
- \frac{1}{\beta} \, \kappa \, \frac{u}{| \nabla q |}
+ \frac{1}{\beta} \, u \, \frac{\nabla q^T \, \nabla^2 q \, \nabla q}{| \nabla q |^4},
\end{eqnarray}
where we set
\begin{equation}
\label{eq:def_u}
u = \nabla \xi \cdot \nabla q.
\end{equation}
Since $\overline{X}_t$ satisfies the constraint $q \left(
  \overline{X}_t \right) = 0$, the dynamics~\eqref{eq:titi} can be
rewritten only in terms of $\xi \left( \overline{X}_t \right) =: z_t$,
in the form
\begin{equation}
\label{eq:path_cv}
dz_t = 
\widetilde{d}_1(z_t) \, dt + \widetilde{d}_2(z_t) \, dt + 
\sqrt{2 \beta^{-1}} \, \, \widetilde{\gamma}(z_t) \, dB_t
\end{equation}
where, for any $z$,
\begin{equation}
\label{eq:def_tilde}
\hspace{-1.5cm}
\widetilde{d}_1(z) = d_1\left(\chi^{-1}(z,0)\right),
\quad
\widetilde{d}_2(z) = d_2\left(\chi^{-1}(z,0)\right),
\quad
\widetilde{\gamma}(z) = | \nabla \xi | \left(\chi^{-1}(z,0)\right).
\end{equation}

\subsection{Effective dynamics associated to~\eqref{eq:X_eps} using
  conditional expectations}

We now follow the strategy that we have outlined in
Section~\ref{sec:closed}. Starting from~\eqref{eq:X_eps}, we first
compute the time derivative of $\xi(X_t^\eps)$ by an It\^o computation,
and next take the conditional expectations of the drift and the diffusion
terms. We hence obtain the effective dynamics~\eqref{eq:y}, where $b_\eps$
and $\sigma_\eps$ (that depend on $\eps$ since the Gibbs measure
$\mu_\eps$ depends on $\eps$) are defined by~\eqref{eq:b}
and~\eqref{eq:sigma} and read  
\begin{eqnarray}
b_\eps(\xib)
&=&
\E_{\mu_\eps} \left[ 
\left(- \nabla V_\eps \cdot \nabla \xi + \beta^{-1} \Delta \xi \right) 
(X) \ | \ \xi(X)=\xib \right]
\nonumber
\\
&=&
\E_{\mu_\eps} \left[ 
\left(- \nabla V_0 \cdot \nabla \xi + \beta^{-1} \Delta \xi \right) (X) 
\ | \ \xi(X)=\xib \right]
\nonumber
\\
& & \quad \quad
- \frac{2}{\eps} \ \E_{\mu_\eps} \left[
\left( q \, \nabla q \cdot \nabla \xi \right) (X)
\ | \ \xi(X)=\xib \right]
\nonumber
\\&=&
\widetilde{d}_1^\eps(\xib)
- \frac{2}{\eps} \ \E_{\mu_\eps} \left[
\left( q \, \nabla q \cdot \nabla \xi \right) (X)
\ | \ \xi(X)=\xib \right]
\label{eq:b_eps}
\\
\sigma_\eps^2(\xib)
&=&\E_{\mu_\eps} \left(|\nabla \xi|^2(X) \ | \ \xi(X)=\xib \right),
\nonumber
\end{eqnarray}
where $\dis{ 
\widetilde{d}_1^\eps(\xib) = 
\E_{\mu_\eps} \left[ 
\left(- \nabla V_0 \cdot \nabla \xi + \beta^{-1} \Delta \xi \right) (X) 
\ | \ \xi(X)=\xib \right]
}$. It is easy to check that, for any $\xib$, we have 
\begin{equation}
\label{eq:limit}
\widetilde{d}_1^\eps(\xib) = \widetilde{d}_1(\xib) + O(\eps)
\quad \text{and} \quad
\sigma_\eps(\xib) = \widetilde{\gamma}(\xib) + O(\eps),
\end{equation}
where $\widetilde{d}_1$ and $\widetilde{\gamma}$ are defined
by~\eqref{eq:def_tilde}. 

\subsection{Sufficient conditions for the pathwise convergence to the
  effective dynamics~\eqref{eq:y}}
\label{sec:pathwise_theory}

Let us establish sufficient conditions under which the
equation~\eqref{eq:path_cv} is equivalent to the effective
dynamics~\eqref{eq:y}, in the limit $\eps \to 0$. 
We hence request that, in
the limit $\eps \to 0$, the dynamics~\eqref{eq:y} and~\eqref{eq:path_cv}
have the same drift and diffusion coefficients. 

We first see that this condition is satisfied for the diffusion
coefficients, in view of~\eqref{eq:limit}: the diffusion coefficient
$\sigma_\eps$ of~\eqref{eq:y} converges to $\widetilde{\gamma}$, the
diffusion coefficient of~\eqref{eq:path_cv}.

We now turn to the drift terms, which is $\widetilde{d}_1 +
\widetilde{d}_2$ in the case of~\eqref{eq:path_cv}, and $b_\eps$ given
by~\eqref{eq:b_eps} for the effective dynamics~\eqref{eq:y}. 
In view of~\eqref{eq:limit}, these drift terms are equal, in the limit $\eps \to
0$, if and only if
\begin{equation}
\label{eq:lim1}
- \widetilde{d}_2(\xib) = \lim_{\eps \to 0} 
\frac{2}{\eps} \ \E_{\mu_\eps} \left[
\left( q \, \nabla q \cdot \nabla \xi \right) (X)
\ | \ \xi(X)=\xib \right].
\end{equation}
In view of~\eqref{eq:d2} and~\eqref{eq:def_tilde}, we have 
\begin{equation}
\label{eq:relation_d2}
\widetilde{d}_2(\xib) = d_2 \left(\chi^{-1}(\xib,0)\right)
= 
d_{2a}\left(\chi^{-1}(\xib,0)\right) + \beta^{-1} d_{2b}\left(\chi^{-1}(\xib,0)\right),
\end{equation}
where $d_{2a}$ and $d_{2b}$ do not depend on $\beta$:
\begin{eqnarray}
\label{eq:d2a}
d_{2a} &=& 
u \, \frac{\nabla q \cdot \nabla V_0}{| \nabla q |^2},
\\
\label{eq:d2b}
d_{2b} &=& 
- \frac{\nabla q \cdot \nabla u}{| \nabla q |^2}
- \kappa \, \frac{u}{| \nabla q |}
+ u \, \frac{\nabla q^T \, \nabla^2 q \, \nabla q}{| \nabla q |^4}.
\end{eqnarray}
On the other hand, we compute, for any $\alpha$,
\begin{eqnarray*}
\E_{\mu_\eps} \left[
\left( q \, \nabla q \cdot \nabla \xi \right) (X)
\ | \ \xi(X)=\alpha \right]
&=&\E_{\mu_\eps} \left[
q(X) \, u(X) \ | \ \xi(X)=\alpha \right]
\\
&=& \int_{\Sigma_\alpha} q \, u \ d \mu_{\eps,\Sigma_\alpha}.
\end{eqnarray*}
A direct computation shows that
\begin{equation}
\label{eq:reste}
\E_{\mu_\eps} \left[
\left( q \, \nabla q \cdot \nabla \xi \right) (X)
\ | \ \xi(X)=\alpha \right]
=
\frac{\eps}{2\beta} \ {\cal E}(\alpha)
+ O(\eps^{3/2}),
\end{equation}
where ${\cal E}(\alpha)$ does not depend on $\beta$ and reads
$$
{\cal E}(\alpha) = \frac{\partial \widetilde{u}}{\partial q}(\alpha,0) 
+
\frac{\widetilde{u}(\alpha,0)}{j(\alpha,0)} \
\frac{\partial j}{\partial q}(\alpha,0)
+ 
\widetilde{u}(\alpha,0) \ \frac{\partial \widetilde{V}_0}{\partial q}(\alpha,0),
$$
where 
\begin{equation}
\label{eq:def_u_tilde}
\widetilde{u}(\xi,q) = u(\chi^{-1}(\xi,q)),
\end{equation} 
$\widetilde{V}_0(\xi,q) = V_0(\chi^{-1}(\xi,q))$, and $j = \det \mbox{jac }
\chi^{-1}$. Hence, \eqref{eq:lim1} reads
\begin{equation}
\label{eq:lim2_pre}
- d_{2a}(\chi^{-1}(\xib,0)) - \frac{1}{\beta} d_{2b}(\chi^{-1}(\xib,0))
=
\frac{1}{\beta} \ {\cal E}(\xib).
\end{equation}
We want to enforce this relation for any 
$\beta$. Since $d_{2a}$, $d_{2b}$ and ${\cal E}$ do not depend on
$\beta$, this yields
\begin{equation}
\label{eq:lim2}
d_{2a}(\chi^{-1}(\xib,0)) = 0
\quad \mbox{and} \quad
- d_{2b}(\chi^{-1}(\xib,0)) = {\cal E}(\xib).
\end{equation}
In view of~\eqref{eq:d2a}, a sufficient condition for the first relation
to hold is 
\begin{equation}
\label{eq:CS1}
\forall \alpha \in \R, \quad u(\chi^{-1}(\xib,0)) = 0,
\end{equation}
where, we recall, $u = \nabla \xi \cdot \nabla q$ and $\chi$ is such
that $\chi(X) = (\xi(X),q(X))$. In what follows, we now assume that
$\xi$ is such 
that~\eqref{eq:CS1} holds. The second relation of~\eqref{eq:lim2} now
reads
\begin{equation}
\label{eq:CS2}
\forall \alpha \in \R, \quad 
\frac{\nabla q \cdot \nabla u}{| \nabla q |^2}(\chi^{-1}(\xib,0))
=
\frac{\partial \widetilde{u}}{\partial q}(\xib,0).
\end{equation}

We have thus proved the following result:

\begin{proposition}
\label{prop:pathwise}
Consider the two-dimensional dynamics~\eqref{eq:X_eps}, and its
one-dimensional limit~\eqref{eq:Xbar}, when $\eps \to 0$. 
On the other hand, consider the one-dimensional effective
dynamics~\eqref{eq:y}, obtained using conditional expectations, and pass
to the limit $\eps \to 0$ in the drift and diffusion coefficients. 

Under the conditions~\eqref{eq:CS1} and~\eqref{eq:CS2} (where $u$, $\chi$
and $\widetilde{u}$ are defined by~\eqref{eq:def_u}, \eqref{eq:def_chi}
and~\eqref{eq:def_u_tilde} respectively), these two
dynamics are the same. In addition, for any $T > 0$, there exists $C >
0$ and $\eps_0 > 0$ such that, for all $\eps \leq \eps_0$, we have
\begin{equation}
\label{eq:estim_pathwise}
\sup_{t \in [0,T]} \E
\left| \xi\left(X_t^\eps\right) - y_t^\eps \right|^2 \leq C \eps,
\end{equation}
where $X_t^\eps$ solves~\eqref{eq:X_eps} and $y_t^\eps$ solves the effective
dynamics~\eqref{eq:y}. 
\end{proposition}

\begin{proof}
We only have to prove the bound~\eqref{eq:estim_pathwise}.
We infer from~\eqref{eq:estim_tony} and assumption {\bf [H1]} that
\begin{equation}
\label{eq:estim_tony_re}
\sup_{t \in [0,T]} \E
\left| \xi \left( X_t^\eps \right) - \xi \left( \overline{X}_t \right)
\right|^2 \leq C \eps.
\end{equation}
The drift and diffusion coefficients of
the effective dynamics on $y_t^\eps$ are $b_\eps$ and $\sigma_\eps$. In
view of~\eqref{eq:b_eps}, \eqref{eq:limit}, \eqref{eq:reste}, 
\eqref{eq:lim2_pre} and~\eqref{eq:relation_d2}, the former satisfies
\begin{eqnarray*}
b_\eps(\alpha) 
&=& 
\widetilde{d}_1^\eps(\alpha) - 
\frac{2}{\eps} \ \E_{\mu_\eps} \left[
\left( q \, \nabla q \cdot \nabla \xi \right) (X)
\ | \ \xi(X)=\xib \right]
\\
&=&
\widetilde{d}_1(\alpha) + O(\eps) - \beta^{-1} {\cal E}(\alpha) + O(\sqrt{\eps})
\\
&=& \widetilde{d}_1(\alpha) + \widetilde{d}_2(\alpha) + O(\sqrt{\eps}).
\end{eqnarray*}
In view of~\eqref{eq:limit}, the latter satisfies
$$
\sigma_\eps(\alpha) = \widetilde{\gamma}(\xib) + O(\eps).
$$
Hence, the difference between, on the one hand, the drift and diffusion
coefficients of the effective dynamics ($b_\eps$ and $\sigma_\eps$)
and, on the other hand, the drift and diffusion coefficients of the
equation~\eqref{eq:path_cv} on $z_t = \xi \left( \overline{X}_t \right)$
(namely $\widetilde{d}_1 + \widetilde{d}_2$ and $\widetilde{\gamma}$),
is of order $O \left( \sqrt{\eps} \right)$. 
We infer from this estimate that, on any bounded time interval, 
$$
\sup_{t \in [0,T]} \E
\left| y_t^\eps - \xi \left( \overline{X}_t \right)
\right|^2 \leq C \eps.
$$
Gathering that estimate with~\eqref{eq:estim_tony_re}
yields~\eqref{eq:estim_pathwise}.
\end{proof}

\medskip

In Sections~\ref{sec:para} and~\ref{sec:num_D2}, we outlined the
condition $\nabla \xi \cdot \nabla q = 0$ as an important condition to
get a good analytical estimate in the sense of [D3], and good numerical
results in terms of residence times.  
If $u = \nabla \xi \cdot \nabla q = 0$, then
conditions~\eqref{eq:CS1} and~\eqref{eq:CS2} are satisfied, and we also
get pathwise convergence (i.e. accuracy in the sense of [D1]), in the simple two-dimensional setting
considered in this section. 

Hence, the same condition $\nabla \xi \cdot \nabla q = 0$ appears,
independently of the estimator ([D3], [D2] or [D1]) that we choose to
measure the accuracy of the effective dynamics.

\subsection{A standard test-case}

Consider the two-dimensional potential energy
\begin{equation}
\label{eq:V_ex}
V_\eps(x,y) = V_0(x,y) + \frac{\Omega^2(x) \, y^2}{\eps},
\quad x \in \rens, \ y \in \rens,
\end{equation}
where $\Omega$ is bounded away from 0 and $V_0$ does not depend on
$\eps$, and the associated overdamped Langevin equation, which defines the
process $X_t^\eps = (x_t^\eps,y_t^\eps)$. The limit dynamics on
$x_t^\eps$ when $\eps \to 0$ 
is well-known in that case (see for instance~\cite{reich00}): it reads
\begin{equation}
\label{eq:eq_reich}
dx_t = - \left( 
\partial_x V_0(x_t,0) + \frac{\Omega'(x_t)}{\beta \Omega(x_t)} 
\right) \, dt + \sqrt{2 \, \beta^{-1}} \, dB_t,
\end{equation}
which is the overdamped Langevin equation associated to the potential 
$V_{\rm eff}(x) = V_0(x,0) + \beta^{-1} \ln \Omega(x)$.

We now wish to recover that result within our approach. 
The potential energy~\eqref{eq:V_ex} is of the form~\eqref{eq:ex_V},
with $q(x,y) = \Omega(x) \,
y$. We wish to choose the reaction coordinate $\xi(x,y) = x$.
Observe then that $u = \nabla \xi \cdot \nabla q = \Omega'(x) \, y \neq
0$. Hence the simple sufficient condition $u=0$ (see
end of Section~\ref{sec:pathwise_theory}) is not satisfied. However, it is easy to see
that the less demanding conditions~\eqref{eq:CS1} and~\eqref{eq:CS2} are
satisfied.  
Hence, in the limit $\eps \to 0$, the effective dynamics~\eqref{eq:y} is
accurate in the sense of pathwise convergence. 

For $\eps > 0$, the effective dynamics reads
\begin{equation}
\label{eq:eff_dyn_part}
d\xi_t = b_\eps(\xi_t) \, dt + \sqrt{2 \, \beta^{-1}} \, dB_t,
\end{equation}
with
$\dis{ 
b_\eps(\alpha) = 
- \E_{\mu_\eps(\alpha,\cdot)} \left( 
\partial_x V_0(\alpha,y) + 2 \frac{\Omega'(\alpha) \Omega(\alpha) y^2}{\eps}
\right)
}$. 
A straightforward computation shows that 
$\dis{ \lim_{\eps \to 0} b_\eps(\alpha)
=
-\partial_x V_0(\alpha,0) - \frac{\Omega'(\alpha)}{\beta \Omega(\alpha)} .
}$
Inserting this relation in~\eqref{eq:eff_dyn_part}, we
recover~\eqref{eq:eq_reich}. 

Hence, taking the limit $\eps \to 0$ in the effective dynamics that we
propose, we recover a well-known result.

\subsection{Numerical results on the example~\eqref{eq:dw}}

\begin{figure}
\centerline{
\input{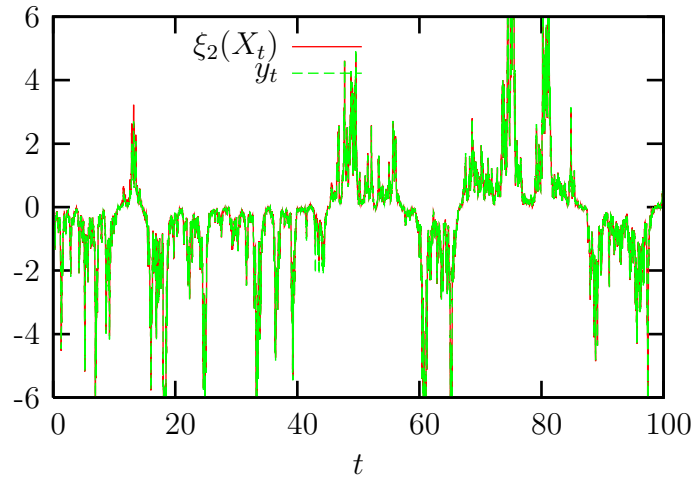}
}
\caption{Comparison of $\xi_2(X_t)$, where $X_t$
  solves~\eqref{eq:X}, and $y_t$ solution of~\eqref{eq:y} with the
  reaction coordinate $\xi_2$.}
\label{fig:traj_xi2}
\end{figure}

In the numerical case considered in Section~\ref{sec:num_D2}, we showed
that the reaction coordinate $\xi_2(x,y) = x \exp(-2y)$ satisfies the
relation $\nabla \xi_2 \cdot \nabla q = 0$. In view of
Proposition~\ref{prop:pathwise}, we hence expect good results when
working with $\xi_2$, in terms of pathwise convergence. We have checked
this as follows. First, we have simulated a solution
of~\eqref{eq:X_eps} (with the potential $V_\eps$ defined
by~\eqref{eq:dw}), for a given realization of the two-dimensional
noise, with $\eps = 0.01$. From this trajectory $X_t$ (we omit
here for clarity the dependence with respect to $\eps$), we
obtain the time evolution $\xi_2(X_t)$, and we can also construct
the one-dimensional noise~\eqref{eq:dB}. This noise is next used
in the effective dynamics~\eqref{eq:y}. We compare both trajectories on
Figure~\ref{fig:traj_xi2}: we observe an excellent agreement over $10^6$
time steps (the trajectories plotted on Figure~\ref{fig:traj_xi2} have
been computed with a time step $\Delta t = 10^{-4}$, hence $T = 100 = 
10^6 \Delta t$).

\begin{figure}
\centerline{
\input{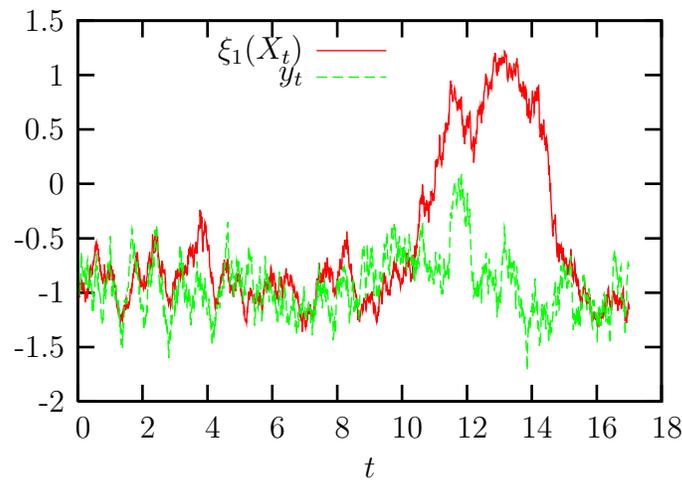}
}
\caption{Comparison of $\xi_1(X_t)$, where $X_t$
  solves~\eqref{eq:X}, and $y_t$ solution of~\eqref{eq:y} with the
  reaction coordinate $\xi_1$.}
\label{fig:traj_xi1}
\end{figure}

In Section~\ref{sec:num_D2}, we also considered the reaction coordinate 
$\xi_1(x,y) = x$, which is such that $u_1(x,y) = \nabla \xi_1 \cdot \nabla q =
2 x \neq 0$. With this choice of reaction coordinate, $\chi^{-1}(\xi,q)
= (\xi,q+1-\xi^2)$, hence $u_1(\chi^{-1}(\xi,0)) = 2 \xi$, so
condition~\eqref{eq:CS1} is not satisfied. We have numerically
performed the same comparison with $\xi_1$ as the one reported above for
$\xi_2$. Results are shown on Figure~\ref{fig:traj_xi1}: we observe that
the complete dynamics (projected on the reaction coordinate) and the
effective dynamics disagree, as expected. Note also the difference in
time ranges between Figures~\ref{fig:traj_xi2} and~\ref{fig:traj_xi1}
(the former corresponding to a time interval 5 times larger than the latter).

\subsection{Numerical results on a three atom molecule}

We conclude this section by considering a system closer to those
considered in molecular simulation, although we acknowledge that it is
still a toy-example. The system is made of three
two-dimensional particles at position $r_i \in \rens^2$, $1 \leq i \leq 
3$ (hence $X=(r_1,r_2,r_3) \in \R^6$), and submitted to the potential
\begin{eqnarray*}
V(X) &=& 
\frac{1}{2\eps} \left( \left\| r_1 - r_2 \right\| - \ell_0 \right)^2
+
\frac{1}{2\eps} \left( \left\| r_2 - r_3 \right\| - \ell_0 \right)^2
+ 
\frac12 k_\theta (\theta(X) - \theta_0)^2
\\
&=&
\frac{1}{2\eps} \left( q_1(X)^2 + q_3(X)^2 \right)
+ 
\frac12 k_\theta (\theta(X) - \theta_0)^2,
\end{eqnarray*}
where $\theta(X)$ is the angle between the bonds $(r_1,r_2)$ and
$(r_2,r_3)$, $q_1(X) = \left\| r_1 - r_2 \right\| - \ell_0$ and $q_3(X) =
\left\| r_2 - r_3 \right\| - \ell_0$. In the above potential,
$\ell_0$ is an equilibrium length whereas $\theta_0$ is an equilibrium
angle.
This potential represents stiff bonds between particles 1 and~2 on
the one hand, and 2 and 3 on the other hand, with a softer term
depending on the three-body angle $\theta$. To remove rigid body motion
invariance, we set $r_2 = 0$ and $r_1 \cdot e_y = 0$. 
Then it is easy to see that the angle $\theta(X)$ satisfies $\nabla
\theta \cdot \nabla q_1 = \nabla \theta \cdot \nabla q_3 = 0$, and hence
seems to be a good reaction coordinate, in view of the several
discussions above. 

\begin{figure}
\centerline{
\input{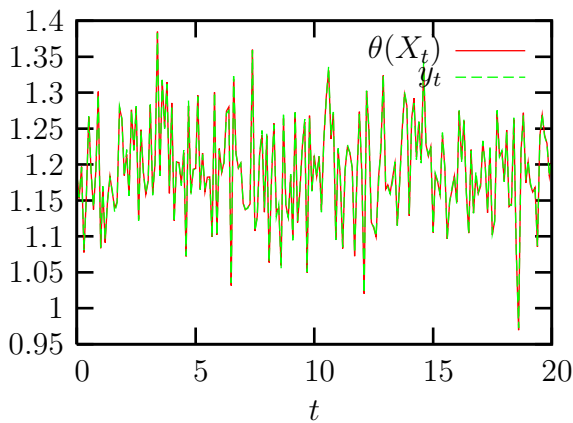}
}
\caption{Comparison of $\theta(X_t)$, where $X_t$
  solves~\eqref{eq:X}, and $y_t$ solution of~\eqref{eq:y} with the
  reaction coordinate $X \mapsto \theta(X)$.}
\label{fig:eau}
\end{figure}

Numerical experiments confirm this belief: choosing this reaction
coordinate, we considered the effective dynamics~\eqref{eq:y}, and
compared its solution with the time evolution $\theta(X_t)$, where $X_t$
solves~\eqref{eq:X}. Results are shown on Figure~\ref{fig:eau} (we
worked with the numerical parameters $\eps = 10^{-3}$, $\ell_0 = 1$,
$\theta_0 = 1.187$ and $k_\theta = 208$): again, we see a good agreement
between both trajectories. 

\ack

This work is supported in part by the INRIA, through the grant 'Action
de Recherche Collaborative' HYBRID, and by the MEGAS non-thematic
program (Agence Nationale de la Recherche, France).
Part of this work was completed while the two authors were visiting the 
Hausdorff Research Institute for Mathematics (HIM, Bonn), whose
hospitality is gratefully acknowledged.
We also wish to acknowledge very enlightning discussions with F. Otto
and G. Menz.

\section*{References}

\bibliographystyle{plain}
\bibliography{legoll_lelievre}

\end{document}